\documentclass[aop,MSNbibl,dvips]{arximspdf}

\doi{10.1214/13-AOP866} 
\volume{42}
\issue{6}
\pubyear{2014}
\firstpage{2454}
\lastpage{2496}

\makeatletter

\newcommand{\rrvert}{\vert}
\newcommand{\llvert}{\vert}
\newtheorem{theorem}{Theorem}[section]
\newproclaim{definition}[theorem]{Definition}
\newtheorem{lemma}[theorem]{Lemma}
\newtheorem{proposition}[theorem]{Proposition}

\newproclaim{remark}[theorem]{Remark}
\newproclaim{example}[theorem]{Example}
\newproclaim{assumption}{Assumption}

\newcommand{\frg}{\mathbf{g}}

\newcommand{\dC}{\mathbb{C}}
\newcommand{\dE}{\mathbb{E}}
\newcommand{\dN}{\mathbb{N}}
\newcommand{\dP}{\mathbb{P}}
\newcommand{\dR}{\mathbb{R}}
\newcommand{\dS}{\mathbb{S}}
\newcommand{\dZ}{\mathbb{Z}}
\newcommand{\cA}{\mathcal{A}}
\newcommand{\cB}{\mathcal{B}}
\newcommand{\cD}{\mathcal{D}}
\newcommand{\cF}{\mathcal{F}}
\newcommand{\cG}{\mathcal{G}}
\newcommand{\cH}{\mathcal{H}}
\newcommand{\cP}{\mathcal{P}}
\newcommand{\cS}{\mathcal{S}}
\newcommand{\cX}{\mathcal{X}}
\renewcommand{\a}{\alpha}
\renewcommand{\b}{\beta}
\renewcommand{\d}{\delta}
\newcommand{\veps}{\varepsilon}
\newcommand{\g}{\gamma}
\renewcommand{\l}{\lambda}
\renewcommand{\o}{\omega}
\renewcommand{\r}{\rho}
\renewcommand{\t}{\tau}
\newcommand{\si}{\sigma}
\newcommand{\grad}{\nabla}
\newcommand{\wt}{\widetilde}
\newcommand{\IND}{\mathbf{1}}
\newcommand{\TR}{\operatorname{tr}}
\newcommand{\supp}{\operatorname{supp}}
\newcommand{\rank}{\operatorname{rank}}
\newcommand{\VAR}{\operatorname{Var}}

\renewcommand{\Im}{\operatorname{\mathfrak{Im}}}
\renewcommand{\Re}{\operatorname{\mathfrak{Re}}}
\renewcommand{\deg}{\operatorname{deg}}
\newcommand{\weak}{\rightsquigarrow}
\newcommand{\locweak}{\stackrel{\mathrm{loc}}{\rightsquigarrow}}
\newcommand{\projweak}{\stackrel{\mathrm{proj}}{\rightsquigarrow}}
\newcommand{\projto}{\stackrel{\mathrm{proj}}{\to}}
\newcommand{\locto}{\stackrel{\mathrm{loc}}{\to}}
\newcommand{\dloc}{d_{\mathrm{loc}}}
\newcommand{\dproj}{d_{\mathrm{proj}}}
\newcommand{\proj}{\mathrm{proj}}
\newcommand{\loc}{\mathrm{loc}}
\newcommand{\sym}{\mathrm{sym}}

\newcommand{\wchi}{\widetilde\chi}
\newcommand{\wcG}{\widetilde{\mathcal{G}}}
\newcommand{\eqref}[1]{(\ref{#1})}
\makeatother

\begin{document}
\begin{frontmatter}

\title{A large deviation principle for Wigner matrices without Gaussian tails\thanksref{T1}}
\runtitle{Large deviations principle for Wigner matrices}
\thankstext{T1}{Supported by the GDRE GREFI-MEFI CNRS-INdAM and supported in part by the European Research
Council through the ``Advanced Grant'' PTRELSS 228032.}

\begin{aug}
\author[A]{\fnms{Charles} \snm{Bordenave}\corref{}\ead[label=e1]{charles.bordenave@math.univ-toulouse.fr}\ead[url,label=u1]{http://www.math.univ-toulouse.fr/\textasciitilde  bordenave}}
\and
\author[B]{\fnms{Pietro} \snm{Caputo}\ead[label=e2]{caputo@mat.uniroma3.it}\ead[url,label=u2]{http://www.mat.uniroma3.it/users/caputo}}
\runauthor{C. Bordenave and P. Caputo}
\affiliation{IMT UMR 5219 CNRS and Universit\'e Paul-Sabatier Toulouse III,\break
and Universit\`a Roma Tre}
\address[A]{Institut de Math\'ematiques de Toulouse\\
CNRS and Universit\'e Toulouse III\\
118 route de Narbonne\\
31062 Toulouse\\
France\\
\printead{e1}\\
\printead{u1}} 
\address[B]{Dipartimento di Matematica e Fisica\\
Universit\`{a} Roma Tre\\
Largo San Murialdo 1\\
00146 Roma\\
Italy\\
\printead{e2}\\
\printead{u2}}
\end{aug}

\received{\smonth{8} \syear{2012}}
\revised{\smonth{3} \syear{2013}}

%
\begin{abstract}
We consider $n\times n$ Hermitian matrices with i.i.d. entries
$X_{ij}$ whose
tail probabilities
$\dP(|X_{ij}|\geq t)$ behave like $e^{-a t^\a}$ for some $a>0$ and
$\a
\in(0,2)$.
We establish a large deviation principle for the empirical spectral
measure of $X/\sqrt{n}$
with speed $n^{1+\a/2}$ with a good rate function $J(\mu)$ that is
finite only if
$\mu$ is of the form $\mu=\mu_{\mathrm{sc}} \boxplus\nu$ for some probability
measure $\nu$ on $\dR$, where
$\boxplus$ denotes the free convolution and $\mu_{\mathrm{sc}}$ is Wigner's
semicircle law. We obtain explicit expressions for
$J(\mu_{\mathrm{sc}} \boxplus\nu)$ in terms of the $\a$th moment of $\nu$.
The proof is based on the analysis of large deviations for the
empirical distribution of
very sparse random rooted networks.
\end{abstract}

%

\end{frontmatter}

\section{Introduction}\label{sec:intro}

Let $\cH_n(\dC)$ denote the set of $n\times n$ Hermitian matrices. The
\emph{empirical spectral measure} of a matrix $A\in\cH_n(\dC)$ is the
probability measure on~$\dR$ defined by
\[
\mu_A = \frac{1} n \sum_{k=1}
^n \delta_{\lambda_k (A)}, %
\]
where $\lambda_1(A) \geq\cdots\geq\lambda_n(A)$ denote the
eigenvalues of $A$ counting multiplicity.
Below, we consider the empirical spectral measure of a \emph{Wigner
random matrix} $X$ described as follows. Let $(X _{ij})_{1 \leq i<j}$
be i.i.d. complex random variables
with variance $\dE|X_{12} - \dE X_{12} |^2 =1$, and let $(X_{ii})_{i
\geq1}$
be 
an independent family of i.i.d. real random variables.
Extend this array by setting
$X_{ij}= \overline X_{ji}$ for $ 1\leq j < i$, and consider the
sequence of
$n\times n$ Hermitian random matrices
%
\begin{equation}
\label{wmatrix} X(n) = (X_{ij}) _{ 1 \leq i, j \leq n}.
\end{equation}
For ease of notation, we often drop the argument $n$ and simply write
$X$ for $X(n)$.

The space $\cP(\dR)$ of probability measures on $\dR$ is endowed with
the topology of weak convergence: a sequence of probability measures
$(\mu_n)_{n \geq1}$ \emph{converges weakly} to $\mu$ if for any
bounded continuous function $f \dvtx \dR\mapsto\dR$,\break
$\int f \,d\mu_n \to\int f \,d\mu$
as $n$ goes to infinity.
We denote this convergence by $\mu_n
\weak \mu$.
Wigner's celebrated theorem asserts that almost surely,
%
\begin{equation}
\label{wsc} \mu_{X / \sqrt n} \weak \mu_{\mathrm{sc}},
\end{equation}
where $\mu_{\mathrm{sc}}$ is the semicircle law, that is, the probability
measure with density $\frac{1}{2\pi}\sqrt{ 4 - x^2}$
on $[-2,2]$; see, for example, \cite{MR2567175,AGZ,MR2808038}.

We consider large deviations, that is, events of the form $\mu_{X /
\sqrt n} \in B$ where $B$ is a measurable set in $\cP(\dR)$ whose
closure does not contain the limiting law $ \mu_{\mathrm{sc}}$.
Clearly, \eqref{wsc} implies that
$\dP(\mu_{X / \sqrt n} \in B)\to0$, $n\to\infty$.
It follows from known concentration estimates that if
the entries $X_{ij}$ are bounded, or if they satisfy a logarithmic
Sobolev inequality, then $\dP(\mu_{X / \sqrt n} \in B)$ decays to $0$
as fast as $e^{-c n^2}$ for some constant $c>0$; see Guionnet and
Zeitouni \cite{MR1781846} or \cite{AGZ}. Further, if the $X_{ij}$ have
a Gaussian law such that $X$ belongs to the Gaussian unitary ensemble
GUE or the Gaussian orthogonal ensemble GOE, then a full large
deviation principle for $\mu_{X / \sqrt n}$ with speed $n^2$ has been
established by Ben Arous and Guionnet in \cite{MR1465640}.
However, apart from the GUE and GOE cases, we are not aware of any case
for which the large deviation principle\vspace*{1pt} for $\mu_{X / \sqrt n}$ has
been obtained. We refer to the recent work of Chatterjee and Varadhan
\cite{MR2890846} for the large deviations of the largest eigenvalues of
$X/n$. For other models of random matrices where the joint law of the
eigenvalues has a tractable form, large deviation principles have been
proved; see, for example, \cite{AGZ}, Section~2.6, or Eichelsbacher,
Sommerauer and Stolz \cite{MR2849043}.

In this paper, we prove a large deviation principle under the
assumption that $X_{ij}$ has tail probabilities $\dP(|X_{ij}|\geq t)$
of order $e^{-at^\a}$ for some $a>0$, and $\a\in(0,2)$. Before stating
our assumptions and results in detail, let us make some preliminary remarks.

It is not hard to see why $n^{1+\a/2}$ is the natural speed for large
deviations in our setting.
For instance, for a fixed $x\in\dR$, consider the event $|X_{ii}|
\sim
x\sqrt n$, for all $i=1,\ldots,n$, which has probability $e^{-cn^{1+\a
/2}}$, for some $c>0$.
This event forces all eigenvalues of $X/\sqrt n$ to shift by $x$ and, therefore,
produces a shift by $x$ of the
limiting spectral measure $\mu_{\mathrm{sc}}$.
Similarly, by considering deviations on the scale $\sqrt n$ of few
elements $X_{ij}$ in each row of the matrix $X$,
one expects to be able to produce more general deformations
of $\mu_{\mathrm{sc}}$ at a cost of order $n^{1+\a/2}$ on the exponential scale.
It turns out that this picture is correct, provided
the deformations of $\mu_{\mathrm{sc}}$ are of the form
$\mu= \mu_{\mathrm{sc}} \boxplus\nu$ for some $\nu\in\cP(\dR)$, where
$\boxplus$ denotes the \emph{free convolution}. Roughly speaking, the
idea is that the entries of $X$ that are visible on a scale $\sqrt n$
form a very sparse weighted random graph or random network $G_n$ that
is asymptotically independent from the rest of the matrix, and a large
deviation principle for
$\mu_{X / \sqrt n}$ can be deduced from a large deviation principle for
the law of the random network $G_n$. This approach also allows us to
obtain explicit expressions for the rate function. The strategy of
proof developed in the present work for Wigner matrices could certainly
be generalized to other models such as random covariance matrices or
random band matrices with the same type of tail assumptions on the
entries. Large deviations with speed $n^ \alpha$ of the largest
eigenvalue may also be handled with similar techniques.

%

\subsection*{Main result}
We recall that a sequence of random variables $(Z
_n)_{n \geq1}$ with values in a topological space $\cX$ with Borel
$\si
$-field $\cB$, satisfies the \emph{large deviation principle} (LDP) with
rate function $J$ and speed $v$, if
$J\dvtx \cX\mapsto[0,\infty]$ is a lower semicontinuous function,
$v\dvtx \dN\mapsto[0,\infty)$ is a function
which increases to infinity, and
for every $B\in\cB$:
%
\begin{eqnarray}
\label{eq:defLDP} -\inf_{x \in B^\circ}J(x )&\leq&\liminf_{n\rightarrow
\infty}
\frac{1}{v(n)} \log\dP ( Z_n \in B ) \leq\limsup
_{n\rightarrow\infty}\frac{1}{v(n)}\log\dP ( Z_n \in B )
\nonumber
\\[-8pt]
\\[-8pt]
\nonumber
&\leq&-
\inf_{ x \in\overline{B}}J(x ),
\end{eqnarray}
where $B^\circ$ denotes the interior of $B$ and
$\overline{B}$ denotes the closure of $B$.~We recall that the lower
semicontinuity of $J$ means that the level sets
$\{x\in\cX\dvtx  J(x)\leq t\}$, $t\in[0,\infty)$,
are closed subsets of $\cX$. When the level sets are compact, the rate
function $J$ is said to be \emph{good}.

We now introduce our statistical assumption. Let $a,\a\in(0,\infty)$.
We say that a complex
random variable $Y$
belongs to the class $\cS_\alpha(a)$, and write $Y\in\cS_\a(a)$,
if
%
\begin{equation}
\label{hypoa} \lim_{t \to\infty} - t^{-\alpha} \log\dP{{\bigl(|Y |
\geq t \bigr)}} = a,
\end{equation}
and if $Y/|Y|$ and $|Y|$ are independent for large values of $|Y|$,
that is, there exists $t_0>0$ and a probability
$\vartheta\in\cP(\dS^1)$ on the unit circle $\dS^1$ such that for all
$t\geq t_0$, all measurable sets $U\subset\dS^1$, one has
%
\begin{equation}
\label{hypoaa} \dP{{\bigl(Y / |Y | \in U \mbox{ and } |Y |\geq t \bigr)}}=\vartheta(U) \dP
{{\bigl(|Y |\geq t \bigr)}}.
\end{equation}
For instance, if $Y$ is Weibull, that is, $Y$ is a nonnegative random
variable with distribution function $F(t)=1-e^{-at^\a}$, with $\a>0$,
and $a>0$, then
$Y\in\cS_\alpha(a)$, with $\vartheta=\d_1$, the unit mass at the
point $1$.
Clearly, if $Y\in\cS_\alpha(a)$ is real valued, then the associated
measure $\vartheta$ must have support in $\{-1,1\}$.
It will be convenient to allow the value $a=\infty$ in \eqref{hypoa}.
Namely, for $\a>0$ we write $Y\in\cS_\alpha(\infty)$ whenever~\eqref
{hypoa} holds with $a=\infty$. 
We do not require \eqref{hypoaa} in this case. For instance, if $Y$ is
a bounded random variable, then $Y\in\cS_\alpha(\infty)$, for all
$\a
>0$, and if $Y$ has a Gaussian tail, then
$Y\in\cS_\alpha(\infty)$, for all $\a\in(0,2)$.
Moreover, if
$Y\in\cS_\alpha(a)$ for some $\a,a>0$, then
$Y\in\cS_\beta(\infty)$ for all $\b\in(0,\a)$.
We remark that \eqref{hypoaa} is a mild technical condition that we do
not expect to be crucial. However, it will turn out to be convenient
for the analysis of random networks in Section~\ref{sec:esgraphs} below.

Throughout the paper, we assume that the array $\{X_{ij}\}$ is given as
above, that is, we have two independent families of random variables:
the off-diagonal
entries $X_{ij}$, $i<j$, which are i.i.d. copies of a complex random
variable $X_{12}$ with unit variance, and the on-diagonal entries
$X_{ii}$, which are i.i.d. copies of a real random variable $X_{11}$.
The matrix $X=X(n)$ is defined as in \eqref{wmatrix}.
Moreover, the following main assumption
will always be understood without explicit mention.

%
\begin{assumption}\label{A1}
There exist $\a\in(0,2)$ and $a,b\in(0,\infty]$ such that
$X_{12}\in\cS
_\alpha(a)$
and $X_{11}\in\cS_\a(b)$.
\end{assumption}

The main result can be formulated as follows.

%
\begin{theorem}
\label{th:main}
Fix $\a\in(0,2)$ as in Assumption~\ref{A1}. The measures $\mu
_{X/\sqrt
n }$ satisfy the LDP with speed $n^{1 + \alpha/2}$ and good rate
function 
%
\begin{equation}
\label{rate1} J( \mu) =\cases{ \Phi( \nu), &\quad $\mbox{if } \mu= \mu_{\mathrm{sc}}
\boxplus\nu\mbox{ for some $\nu\in \cP (\dR)$},$ \vspace*{2pt}
\cr
\infty,&\quad $ \mbox{otherwise}$,}
\end{equation}
where $\Phi\dvtx \cP(\dR)\mapsto[0,\infty]$ is a good rate function.
\end{theorem}
%
More details on the rate function $\Phi$ will be given in Theorems~\ref
{th:rate} and~\ref{th:rateb} below. We anticipate that $\Phi
(\nu
)=0$ if and only if $\nu= \d_0$, where $\d_0$ is the Dirac mass at $0$.
Moreover, as one should expect, in the case $a=b=\infty$, one has
$\Phi
(\nu)=\infty$ for all $\nu\neq\d_0$.

The proof of Theorem~\ref{th:main} consists of two main parts.
The first part, the ``random matrix theory part'' of the work, is
discussed in Section~\ref{sec:expeq}.
Here, we show that at speed $n^{1+\alpha/2}$ the large deviations are
governed by the sparse $n\times n$ random matrix $C=C(n)$ defined by
\[
C_{ij}=\cases{\displaystyle\frac{X_{ij}}{\sqrt n}, & \quad$\mbox{if }\displaystyle \veps(n)\leq
\frac{X_{ij}}{\sqrt n} \leq\veps (n)^{-1},$\vspace*{2pt}
\cr
0, & \quad$\mbox{otherwise},$ }
\]
where $\veps(n)$ is a cutoff sequence that for convenience will be set
equal to $1/\log n$. In particular, we show that as far as the LDP with
speed $n^{1+\alpha/2}$ is concerned, $\mu_{X/\sqrt n }$
behaves as $\mu_{\mathrm{sc}}\boxplus\mu_{C}$, where $ \mu_{C}$ is the spectral
measure of the matrix
$C$; see Proposition~\ref{prop:exptight} below. As a consequence, the
LDP for $\mu_{X/\sqrt n }$ will be obtained by contraction if one has
the LDP for $\mu_C$ with speed $n^{1+\alpha/2}$ and rate function
$\Phi$.

The second part, the ``random graph theory part'' of the work,
is presented in Section~\ref{sec:esgraphs}. Here, we prove the above
mentioned LDP for the spectral measures~$\mu_C$.
By viewing the matrix $C$ as the adjacency matrix of a weighted graph,
one runs naturally into
the analysis of large deviations for
sparse random networks. This is best formulated within the theory
of local convergence for networks that was recently developed by
Benjamini and Schramm \cite{bensch}, Aldous and Steele \cite
{aldoussteele} and Aldous and Lyons
\cite{aldouslyons}.
Let us briefly sketch the main ideas---all details will be given in
Section~\ref{sec:esgraphs}. Let $G_n$ be the sparse random network
naturally associated to the $n\times n$ matrix $C$, that is, $G_n$ is
the weighted graph with $n$ vertices whose adjacency matrix is given by
$C$. Notice that the weights can have a sign, and there are loops
corresponding to nonzero diagonal entries of $C$. Take a vertex at
random, call it the root, and consider the
connected component of $G_n$ at that vertex.
This gives rise to a random connected rooted network, we call $\r_n$
its law.
By identifying two networks which differ only by a
permutation of the vertex labels,
the law $\r_n$ is regarded as an element of the space $\cP(\cG_*)$
of probability measures on $\cG_*$, where $\cG_*$ is
the space of equivalence classes (under rooted isomorphisms) of
connected rooted networks.
The essential point is that the eigenvalue distribution $\mu_C$ can be
identified with a
suitable ``spectral measure'' $\mu_{\r_n}$ associated to the law $\r
_n$; see also~\cite{MR2724665,MR2789584,MR2857250}
for recent works based on the same idea.

Since the network $G_n$ is very sparse, one has that almost surely $\r
_n$ converges (under the weak local convergence \cite{aldouslyons}) to
the Dirac mass
on the trivial element of $\cG_*$, namely the network consisting of a
single isolated vertex (the root).
We introduce a suitable weak topology on $\cP(\cG_*)$, and prove that
the measures
$\r_n$ satisfy a LDP
with speed $n^{1+\alpha/2}$ and a good rate function $I(\r)$. The
latter is finite only if $\r$ belongs to the so called sofic measures,
that is, if $\r$ is the weak local limit of finite networks, and if
the support of $\r$ satisfies some natural constraints. Call $\cP
_s(\cG
_*)$ the set of such probability measures. We find that for $\r\in\cP
_s(\cG_*)$, one has
%
\begin{equation}
\label{rate2} I( \r) = b \dE_\rho\bigl| \omega_G ( o )
\bigr|^\alpha+ \frac{a}{2} \dE_\rho\sum
_{v \in V_G \setminus o } \bigl|\omega_G ( o,v) \bigr|^\alpha,
\end{equation}
where $\dE_\rho$ denotes expectation w.r.t. $\r$, the law of the
equivalence class of a connected rooted network $(G,o)$,
$o$ denoting the root; $\omega_G ( o )$ denotes the weight of the loop
at the root, and $\omega_G ( o,v)$ denotes the weight of
the edge $(o,v)$ if $v$ is an element of the vertex set $V_G$ of the
network. We refer to Proposition~\ref{prop:LDPnetwork} for the precise result.

It turns out that the choice of a ``myopic'' topology on $\cP(\cG_*)$
is crucial to have the desired result. On the other hand, we want this
topology to be fine enough to have that the map $\r\mapsto\mu_\r$
defining the spectral measure associated to $\r$ is continuous.
If all this is satisfied, then a LDP for the spectral measure $\mu
_C=\mu
_{\r_n}$ can be obtained by contraction from the LDP for $\r_n$; see
Proposition~\ref{ldpmuc}.
In particular, we find that the function $\Phi$ in Theorem~\ref
{th:main} is given by
%
\begin{equation}
\label{rate3} \Phi(\nu)=\inf\bigl\{I(\r), \r\in\cP_s(\cG_*)\dvtx
\mu_\r=\nu\bigr\}.
\end{equation}

We now turn to more explicit characterizations of the rate function in
Theorem~\ref{th:main}. From the approach discussed above, we will see
that the rate function $\Phi$ depends on the laws of $X_{11}$ and
$X_{12}$ only through $\alpha, a, b$ and the supports of the
associated measures on $\dS^1$. While the variational principle \eqref
{rate3} is not always explicitly solvable,
there is a large class of $\nu\in\cP(\dR)$ for which $\Phi(\nu)$
can be
computed. This allows us to give explicit expressions for the rate
function $J(\mu)$ in Theorem~\ref{th:main}.
Recall that the free convolution with $\mu_{\mathrm{sc}}$ is injective: for any
$\mu\in\cP(\dR)$ there is at most one $\nu\in\cP(\dR)$ such that
$\mu= \mu_{\mathrm{sc}} \boxplus\nu$. Let $\cP_\sym(\dR)$ denote the set of
symmetric probability measures on $\dR$. If $\mu= \mu_{\mathrm{sc}} \boxplus
\nu$, then $\mu\in\cP_\sym(\dR)$ is equivalent to $\nu\in\cP
_\sym
(\dR)$.
For more details on free convolution with the semicircular
distribution, we refer to Biane \cite{MR1488333}.
For $\nu\in\cP(\dR)$, we use the notation
%
\begin{equation}
\label{ath} m_\a(\nu)= \int|x|^\alpha\,d \nu(x)
\end{equation}
for the $\a$th moment of $\nu$.
If $X_{11}\in\cS_\a(b)$ for some $b<\infty$, then we write
$\vartheta
_b$ for the associated measure given in \eqref{hypoaa}. 
Recall that since $X_{11}$ is real, $\vartheta_b$ is a measure on $\{
-1,1\}$.
The following theorem summarizes the main facts we can establish about
the rate function.

%
\begin{theorem}\label{th:rate}
\textup{(a)} For any $\nu\in\cP(\dR)$,
\[
\Phi(\nu)\geq{{\biggl(\frac{a} 2 \wedge b \biggr)}} m_\a(
\nu).\vspace*{-6pt}
\]
\begin{longlist}[(b)]
\item[(b)]
If $b<\infty$ and $\supp(\vartheta_b ) = \{ -1, 1\} $,
then for any $\nu\in\cP(\dR)$:
\[
\Phi(\nu)\leq b m_\a(\nu).
\]

\item[(c)] If $b<\infty$ and $\supp(\vartheta_b ) = \{-1,1 \}$, and $\nu
\in
\cP_\sym(\dR)$, then
\[
\Phi(\nu)= {{\biggl(\frac{a} 2 \wedge b \biggr)}} m_\a(\nu).
\]
\end{longlist}
\end{theorem}
Some remarks about Theorem~\ref{th:rate}.
Part (a) shows clearly that $\Phi(\nu) = 0$ is equivalent to $\nu= \d
_0$, that is, $J(\mu) = 0$ is equivalent to $\mu= \mu_{\mathrm{sc}}$. It also
shows that $J$ is a good rate function since the level sets $\{m_\a
(\cdot)\leq t\}$, $t\in[0,\infty)$ are compact in $\cP(\dR)$.
Concerning the remaining statements, the fact that the moments $m_\a
(\nu
)$ appear naturally in the rate function and the special role played by
symmetric measures $\nu$ can be understood as follows.

As one could expect, there is a natural way to achieve a large deviation
$\mu_{X/\sqrt n}\sim\mu_{\mathrm{sc}}\boxplus\nu$ by tilting only the diagonal
entries of $X$, namely by considering events of the form $\mu_{\cD
/\sqrt n}\sim\nu$, where
$\cD$ denotes the diagonal matrix with entries $X_{11},\ldots
,X_{nn}$, and
\[
\mu_{{\cD}/{\sqrt n}}=\frac{1}n\sum_{i=1}^n
\d_{
{X_{i,i}}/{\sqrt n}}.
\]
In view of \eqref{hypoaa}, one can consider an arbitrary $\nu\in\cP
(\dR
)$ here if $b<\infty$ and $\supp(\vartheta_b ) = \{-1,1 \}$. If
$b<\infty$ and $\supp(\vartheta_b ) = \{+1 \}$ (or $\{-1\}$) then only
$\nu$ whose support is $\dR_+$ (or $\dR_-$) can be considered. If
$b=\infty$, then no measure $\nu\neq\d_0$ will have a finite cost on
the scale $n^{1+\a/2}$.

Similarly, one can try to reach a large deviation
$\mu_{X/\sqrt n}\sim\mu_{\mathrm{sc}}\boxplus\nu$ by tilting only the
off-diagonal entries of $X$. For instance,
for $n$ even, let $\cA$ denote the
block diagonal matrix made up of the $2\times2$ blocks
\[
\left(\matrix{ 0 & X_{i,i+1}\vspace*{2pt}
\cr
\bar
X_{i,i+1} & 0} \right),\qquad
i=1,\ldots,n/2.
\]
That is, $\cA$ is defined
by $\cA_{2i-1,2i}=X_{i,i+1}$, $\cA_{2i,2i-1}=\bar X_{i,i+1}$,
$i=1,\ldots
,n/2$, and $\cA_{i,j}=0$ for all other entries.
It is straightforward to see that
the empirical spectral measures 
of $\cA/\sqrt n$ is given by
\[
\mu_{{\cA}/{\sqrt n}}=\frac{1}n\sum
_{i=1}^{n/2} (\d_{
{|X_{i,i+1}|}/{\sqrt n}}+\d_{-{|X_{i,i+1}|}/{\sqrt n}} ).
\]
Notice that $\mu_{\cA/\sqrt n}$ is a symmetric distribution. Thus, if
we try to obtain $\mu_{X/\sqrt n}\sim\mu_{\mathrm{sc}}\boxplus\nu$ by
requiring $\mu_{\cA/\sqrt n}\sim\nu$ we are forced to restrict to
$\nu
\in\cP_\sym(\dR)$.

In view of this discussion, it is natural to look for upper bounds on
the rate function $\Phi$ in terms of the rate function associated to
large deviations of $\mu_{\cD/\sqrt n}$ and $\mu_{\cA/\sqrt n}$.
Our results will show in particular that if the variables $X_{ij}$ are
as in Assumption~\ref{A1}, with $b<\infty$ and $\supp(\vartheta_b ) =
\{-1,1 \}$, then:

\begin{longlist}[(1)]
\item[(1)]$\mu_{\cD/\sqrt n}$ satisfies a LDP on $\cP(\dR)$ with speed
$n^{1+\a/2}$ and rate function $I_b(\nu)=b m_\a(\nu)$,
for all $\nu\in\cP(\dR)$;

\item[(2)]
$\mu_{\cA/\sqrt n}$
satisfies a LDP on $\cP(\dR)$ with speed
$n^{1+\a/2}$ and rate function equal to $I_a(\nu)=\frac{a}2 m_\a
(\nu)$,
for all $\nu\in\cP_\sym(\dR)$, and $I_a(\nu)=+\infty$ if $\nu
\notin\cP
_\sym(\dR)$.
\end{longlist}

Since $\mu_{\cD/\sqrt n}$ and $\mu_{\cA/\sqrt n}$ are the empirical
measures induced by i.i.d. random variables rescaled by $\sqrt n$, the
statements above can be seen as extremal instances of Sanov's theorem,
in the case of variables with exponential tails of the form \eqref{hypoa}.
Thus, roughly speaking, part (b) in Theorem~\ref{th:rate}
can be interpreted as the bound obtained by adopting the strategy $\mu
_{\cD/\sqrt n}\sim\nu$ to reach the deviation $\mu_{X/\sqrt n}\sim
\mu
_{\mathrm{sc}}\boxplus\nu$.
When $b\leq a/2$, parts (a) and~(b) above yield the expression
\[
\Phi(\nu)= b m_\a(\nu),
\]
showing that this strategy is optimal.
Similarly, to illustrate part (c), observe that
if $\nu\in\cP_\sym(\dR)$, then for
the deviation
$\mu_{X/\sqrt n}\sim\mu_{\mathrm{sc}} \boxplus\nu$
one can also use the strategy $\mu_{\cA/\sqrt n}\sim\nu$.
This reasoning will produce the bound
$\Phi(\nu)\leq(a/ 2 \wedge b ) m_\a(\nu)$. The general bound in
part (a)
then shows that this is actually an optimal strategy if $a/2\leq b$.

If the support of $\vartheta_b$ is only $\{+1\}$ (or $\{-1\}$) then
the above scenario changes in that one can use the diagonal matrix $\cD
$ only to reach deviations $\nu$ whose support is $\dR_+$ (or $\dR_-$).
In this case, we have the following estimates. Without loss of
generality, we restrict to $\supp(\vartheta_b ) = \{+1 \}$.

\begin{theorem}\label{th:rateb}
Suppose $b<\infty$, and
$\supp(\vartheta_b) = \{+1\} $.

\begin{longlist}[(a)]
\item[(a)]
If $\supp(\nu) \subset\dR_{+}$,
then
\[
\Phi(\nu)\leq b m_\a(\nu).
\]

\item[(b)] Suppose $\a\in(1,2)$. If $\nu\in\cP_\sym(\dR)$, then
\[
\Phi(\nu)= \frac{a} 2 m_\a(\nu).
\]

\item[(c)] Suppose $\a\in(1,2)$. If $\int x\,d\nu(x)<0$ then
$\Phi(\nu)= +\infty$.
\end{longlist}
\end{theorem}
The above result can be interpreted as before by appealing to the large
deviations of $\mu_{\cD/\sqrt n}$ and $\mu_{\cA/\sqrt n}$. In
particular, part (b) shows that since one cannot realize a symmetric deviation
$\nu\in\cP_\sym(\dR)$
using the matrix $\cD$ only, it is less costly to realize it using the
matrix $\cA$ only. Similarly, in part (c), one has that neither $\cD$
nor~$\cA$, nor any other matrix with vanishing trace, can be used to
produce a measure $\nu$ with $\int x \,d\nu(x)<0$ and, therefore, the
rate function must be $+\infty$. We believe that results in parts (b)
and (c) above should hold without the additional condition $\a\in(1,2)$.

The proofs of Theorems \ref{th:rate} and \ref{th:rateb} are given in
Section~\ref{subsec:Th23}.

\section{Exponential equivalences}
\label{sec:expeq}

Throughout the rest of the paper, we fix the cutoff sequence $\veps
(n)$ as
%
\begin{equation}
\label{eq:seqveps} \veps(n) = \frac{1} { \log n }.
\end{equation}
%
We decompose the matrix $X$ as
%
\begin{equation}
\label{decompo} \frac{X}{\sqrt n}= A+B+C+D,
\end{equation}
where the matrices $A,B,C,D$ are defined by
\begin{eqnarray*}
 A_{ij}& =& \IND_{ |X_{i j } |< (\log n)^{2 / \alpha} } \frac{X_{ij}
}{\sqrt n},\qquad
B_{ij} = \IND_{ (\log n)^{2 / \alpha} \leq|X_{i j
} |\leq\veps(n) n^{1/2} } \frac{X_{ij} }{\sqrt n},
\\
 C_{ij}& =& \IND_{ \veps(n) n^{1/2} < |X_{i j } | < \veps(n) ^{-1}
n^{1/2} } \frac{X_{ij} }{\sqrt n},\qquad
D_{ij} = \IND_{\veps
(n)^{-1} n^{1/2} < |X_{i j } | } \frac{X_{ij} }{\sqrt n}.
\end{eqnarray*}
%
The matrix $A$ represents the bulk of the original matrix, while the
matrix $C$ yields the elements that are visible on the scale $\sqrt n$.
The starting point of our analysis (see Lemmas \ref{le:very} and~\ref{le:moder} below) is to show that
the contribution of both $B$ and $D$ is negligible for large deviations
with speed $n^{1+\a/2}$.

We define the distance on $\cP(\dR)$ as
%
\begin{equation}
\label{eq:defdist} d (\mu, \nu) = \sup\bigl\{ \bigl| g_\mu(z ) -
g_{\nu} (z) \bigr| \dvtx \Im(z) \geq2 \bigr\},
\end{equation}
where $g_\mu$ is the \emph{Cauchy--Stieltjes transform} of $\mu$, that
is, for $z \in\dC_+ = \{ z \in\dC\dvtx  \Im(z) > 0 \}$,
%
\begin{equation}
\label{gstie} g_\mu( z) = \int\frac{ \mu(dx) }{x - z}.
\end{equation}
Recall that this distance is a metric for the weak convergence; see,
for example,~\cite{AGZ}, Theorem~2.4.4.
Let also $d_{\mathrm{KS}}$ denote the Kolmogorov--Smirnov distance and let $W_p$
denote the $L^p$-Wasserstein distance; see Appendix \ref{app} below for
the relevant definitions.
The introduction of the distance $d_{\mathrm{KS}}$ is mainly due to the use of
the rank inequality of Lemma~\ref{le:rank}.
The Wasserstein distance on the other hand can be controlled in terms
of the matrix elements thanks to the Hoffman--Wielandt inequality in
Lemma~\ref{le:HW}. We shall relate these distances to the distance
\eqref{eq:defdist} via the following estimate, which is a consequence
of \eqref{eq:KSdual} and \eqref{eq:KRdual}:
%
\begin{equation}
\label{eq:boundd} d (\mu, \nu) \leq d_{\mathrm{KS}} ( \mu, \nu) \wedge
W_1 ( \mu, \nu).
\end{equation}

The following proposition is the first major step on the way to prove
Theorem~\ref{th:main}.

\begin{proposition}
\label{prop:exptight}
The random probability measures $\mu_{\mathrm{sc}} \boxplus\mu_{C}$ and $ \mu
_{X / \sqrt n } $ are exponentially equivalent: for any $\delta> 0$,
\[
\limsup_{n \to\infty} \frac{1}{n^{1 + \alpha/ 2}} \log\dP \bigl(d( \mu
_{X / \sqrt n }, \mu_{\mathrm{sc}} \boxplus\mu_{C} ) \geq\delta
\bigr) = - \infty.
\]
\end{proposition}

The rest of this section is devoted to
the proof of Proposition~\ref{prop:exptight}. The strategy is as
follows: we start by showing that the contribution of $D$ in \eqref
{decompo} can be neglected (Lemma~\ref{le:very}), then we show that $B$
can also be neglected (Lemma~\ref{le:moder}). The last step will then
consist in proving that $\mu_{A + C}$ and
$\mu_{\mathrm{sc}} \boxplus\mu_{C}$ are exponentially equivalent.
We note that the assumption \eqref{hypoaa}
is not needed for the proof of Proposition~\ref{prop:exptight}.
Actually, a careful look at the proof shows that
it is sufficient to replace condition~\eqref{hypoa} by the weaker
assumption $\limsup_{t\to\infty}t^{-\alpha} \log\dP{{(|Y |\geq t
)}}<0$; 
see Remark~\ref{weaker} below.

\subsection{Preliminary estimates}
%
\begin{lemma}[(Very large entries)]\label{le:very}
The random probability measures $\mu_{A + B+ C}$ and $ \mu_{X / \sqrt n
} $ are exponentially equivalent: for any $\delta> 0$,
\[
\limsup_{n \to\infty} \frac{1}{n^{1 + \alpha/ 2}} \log\dP{{\bigl(d(
\mu_{X / \sqrt n }, \mu_{A + B + C} ) \geq\delta\bigr)}} = - \infty.
\]
\end{lemma}

\begin{pf}
From \eqref{eq:boundd}, it is sufficient to prove that for any $\delta
> 0$,
\[
\limsup_{n \to\infty} \frac{1}{n^{1 + \alpha/ 2}} \log\dP {{
\bigl(d_{\mathrm{KS}} ( \mu_{X / \sqrt n }, \mu_{A + B + C} ) \geq\delta
\bigr)}} = - \infty.
\]
Then, using the rank inequality Lemma~\ref{le:rank}, it is sufficient
to prove that for any $\d>0$
\[
\limsup_{n \to\infty} \frac{1}{n^{1 + \alpha/ 2}} \log\dP {{\bigl(\rank( D )
\geq\delta n \bigr)}} = - \infty.
\]
However, since the rank is bounded by the number of nonzeros entries of
a matrix, one has
\[
\dP{{\bigl(\rank( D ) \geq2 \delta n \bigr)}} \leq\dP \biggl( \sum
_{ 1 \leq i
\leq j \leq n } \IND \bigl(|X_{ij}| \geq
\veps(n)^{-1} n^{1/2} \bigr) \geq\delta n \biggr).
\]
The Bernoulli variables $\IND( |X_{ij}| \geq\veps(n)^{-1} n^{1/2} ), 1
\leq i \leq j \leq n$, are independent. Also, by assumption \eqref
{hypoa}, their mean value $p_{ij} =\dP(|X_{ij}| \geq\veps(n)^{-1}
n^{1/2} )$ satisfies
\[
p_{ij} \leq p(n): = e^{ - c \veps(n)^{-\alpha} n^{\alpha/2} }
\]
for some $ c > 0$. For our choice of $\veps(n)$ in \eqref{eq:seqveps},
one has $ p(n) = o (1/n^2)$.
Hence, it is sufficient to prove that for any $\d>0$:
\[
\limsup_{n \to\infty} \frac{1}{n^{1 + \alpha/ 2}} \log\dP \biggl( \sum
_{ 1 \leq i \leq j \leq n } {{\bigl(\IND \bigl( |X_{ij}| \geq
\veps(n)^{-1} n^{1/2} \bigr) - p_{ij} \bigr)}} \geq
\delta n \biggr) = - \infty.
\]
Recall Bennett's inequality \cite{bennett}: if $W_i$, $i=1,\ldots,m$ are
independent Bernoulli($p_i$) variables, and $h(x) = (x+1)\log( x+1)
-x$, then one has
%
\begin{equation}
\label{Bennett} \dP \Biggl( \sum_{ i = 1}^{ m }
{{(W_{i} -p_i)}} \geq t \Biggr) \leq \exp \biggl( -
\si^2 h \biggl( \frac{t}{\si^2} \biggr) \biggr)
\end{equation}
with $
\si^2 = \sum_{i = 1} ^ {m} p_i (1 - p_i)$.
In our case, for all $n$ large enough,
\[
\si^2 = \sum_{1 \leq i \leq j \leq n} p_{ij} (1
- p_{ij} ) \leq\frac{
n(n+1) p(n)}{ 2}.
\]
Therefore, using $h(x) \sim x \log x$ as $x \to\infty$,
\begin{eqnarray*}
\dP \biggl( \sum_{ 1 \leq i \leq j \leq n } \bigl(\IND \bigl(
|X_{ij}| \geq\veps(n)^{-1} n^{1/2} \bigr) -
p_{ij} \bigr) \geq\delta n \biggr) &\leq&\exp \biggl(-
\si^2 h \biggl( \frac{ n
\delta} { \si^2 } \biggr) \biggr)
\\
& \leq&\exp \bigl( c_0 n \log{{\bigl({ n p(n) } \bigr)}} \bigr)
\end{eqnarray*}
for some constant $c_0 > 0$ depending on $\delta$. Now, since $n\leq
p(n)^{-1/2}$ for $n$ large, we find that for some $c_1>0$, for all $n$
large enough the last expression is upper bounded by
\[
\exp \biggl( \frac{1}2c_0 n \log p(n) \biggr) \leq\exp
\bigl( - c_1 n^{1 +
\alpha/ 2} \veps(n)^{-\alpha} \bigr).
\]
This proves the claim.
\end{pf}

We now show that the contribution of $B$ in \eqref{decompo} is also negligible.
While Lemma~\ref{le:very} would work for any $\a>0$,
the next results use the fact that $\a\in(0,2)$.

\begin{lemma}[(Moderately large entries)]\label{le:moder}
The random probability measures $\mu_{A + C}$ and $ \mu_{X / \sqrt n }
$ are exponentially equivalent: for any $\delta> 0$,
\[
\limsup_{n \to\infty} \frac{1}{n^{1 + \alpha/ 2}} \log\dP{{\bigl(d(
\mu_{X / \sqrt n }, \mu_{A + C} ) \geq\delta\bigr)}} = - \infty.
\]
\end{lemma}

\begin{pf}
From \eqref{eq:boundd}, Lemma~\ref{le:very} and the triangle
inequality, it is sufficient to check that for any $\delta> 0$,
\[
\limsup_{n \to\infty} \frac{1}{n^{1 + \alpha/ 2}} \log\dP{{
\bigl(W_2 ( \mu_{A + B +C}, \mu_{A + C} ) \geq\delta
\bigr)}} = - \infty,
\]
where $W_2 \geq W_1$ is the $L^2$-Wasserstein distance defined by
\eqref
{eq:defWp}. From the Hoffman--Wielandt inequality Lemma~\ref{le:HW}, it
is sufficient to prove that for any $\delta> 0$,
\[
\limsup_{n \to\infty} \frac{1}{n^{1 + \alpha/ 2}} \log\dP {{\biggl(
\frac { 1} { n} \TR\bigl( B^2 \bigr) \geq\delta\biggr)}} = -
\infty.
\]
We write
\[
\frac{ 1} { n} \TR\bigl( B^2 \bigr) \leq\frac{2} { n^2} \sum
_{ 1 \leq i \leq j
\leq n } |X_{ij}|^2 \IND\bigl( (
\log n)^{2 / \alpha} \leq|X_{i j } |\leq \veps(n) n^{1/2}
\bigr).
\]
Thus, from Markov's inequality, for any $\lambda> 0$,
\[
\dP \biggl( \frac{ 1} { n} \TR\bigl( B^2 \bigr) \geq2 \delta
\biggr)  \leq e^{
-\lambda\delta} \prod_{ 1 \leq i, j \leq n} \dE{{
\bigl[e^{ n^{-2}
\lambda|X_{ij}|^2 \IND( (\log n)^{2 / \alpha} \leq|X_{i j } |\leq
\veps(n) n^{1/2})}\bigr]}}.
\]
To estimate the last expectation, we use the
integration by part formula, for $\mu\in\cP(\dR)$ and $g \in C^1$,
%
\begin{eqnarray}
\label{eq:upIbP} \int_a^b g(x) \,d\mu(x) &=& g(a)
\mu\bigl([a, \infty) \bigr) - g(b) \mu\bigl( (b, \infty)\bigr)
\nonumber
\\[-8pt]
\\[-8pt]
\nonumber
&&{} + \int
_a ^ b g'(x) \mu\bigl( [ x,
\infty) \bigr) \,dx.
\end{eqnarray}
Define the
function
%
\begin{equation}
\label{eq:deff} f(x) = n^{-2}\lambda x^2 - c
x^\alpha.
\end{equation}
Let $\mu$ denote the law of $|X_{ij}|$, and $g(x)=e^{n^{-2}\l x^2}$.
By Assumption~\ref{A1}, there exists a constant $c > 0$ such that
%
\begin{equation}
\label{conx}\mu\bigl([t, \infty) \bigr)=\dP\bigl( |X_{ij} | \geq t \bigr) \leq
\exp\bigl( - c t^\alpha\bigr)
\end{equation}
for all $t$ large enough. In particular, $g(t)\mu([t, \infty) )\leq
e^{ f(t)} $.
From \eqref{eq:upIbP}, it follows that
%
\begin{eqnarray}
\label{eq:uptrB2}&& \dE{{\bigl[e^{ n^{-2} \lambda|X_{ij}|^2 \IND( (\log n)^
{2 / \alpha}
\leq|X_{i j } |\leq\veps(n) n^{1/2})}\bigr]}} \nonumber\\
&&\qquad\leq 1+ \int
_{(\log n)^{2 / \alpha} } ^ {\veps(n) n^{1/2}} g(x) \,d\mu (x)
\nonumber
\\[-8pt]
\\[-8pt]
\nonumber
&&\qquad \leq 1+e^{ f((\log n)^{2 / \alpha}) } +\int_{(\log
n)^{2 / \alpha} } ^ {\veps(n) n^{1/2}}
\frac{ 2 \lambda x }{n^{2}} e^{
f(x) } \,dx
\\
& &\qquad\leq 1+e^{ f((\log n)^{2 / \alpha}) } + \frac{ \lambda\veps(n)^2 }{n} \max_{ x \in[(\log n)^{2 / \alpha}, \veps(n) n^{1/2}]}
e^{ f(x)}.\nonumber
\end{eqnarray}
We choose $\lambda= \frac{1}2c \veps(n)^{\alpha-2}n^{1 + \alpha/
2}$, with the constant $c>0$ given in \eqref{conx}.
Simple computations show that $f(x)$
reaches its maximum for
$x\in[(\log n)^{2 / \alpha},\break \veps(n) n^{1/2}]$ at $x = (\log n)^{2 /
\alpha} $, where it is equal to
\[
\frac{1}2c \veps(n)^{\alpha-2} n^{\alpha/ 2 -1 } (\log n)
^{4 /
\alpha} - c ( \log n) ^2.
\]
Using \eqref{eq:seqveps}, for $n\geq n_0$ this is smaller than $-\frac
{c}2(\log n)^2$.
Therefore, using $1+x\leq e^x$, $x\geq0$, one has that \eqref
{eq:uptrB2} is bounded by $\exp{ [e^ {-({c}/4)(\log n)^2} ]}$
for $n$ large enough.
It follows that
\[
\frac{1}{n^{1 + \alpha/ 2}} \log\dP \biggl(\frac{ 1} { n} \TR\bigl( B^2
\bigr) \geq2 \delta \biggr) \leq- \frac{1}2c \delta
\veps(n)^{\alpha-2} + n^{1-\a/2}e^ {-({c}/4)(\log n)^2}.
\]
The desired conclusion follows.
\end{pf}

For $s > 0$, we define the compact set for the weak topology
\[
K_s = \biggl\{ \mu\in\cP(\dR) \dvtx \int x^2 \,d\mu\leq
s \biggr\}.
\]
For a suitable choice of $s$, we now check that $\mu_{C}$ is in $K_s$
with large probability.

\begin{lemma}[(Exponential tightness estimates)]\label{le:exptight}
\[
\limsup_{n \to\infty} \frac{1}{n^{1 + \alpha/ 2}} \log\dP{{(\mu
_{C} \notin K_{(\log n)^2} )}} = - \infty.
\]
Moreover, if $I = \{ (i,j) \dvtx |X_{ij} | > (\log n)^{2/\alpha} \}$, for
any $\delta> 0$,
\[
\lim_{n \to\infty} \frac{1}{n^{1 + \alpha/ 2}} \log\dP \bigl( | I | \geq\delta
n^{1 + \alpha/2} \bigr) = - \infty.
\]
\end{lemma}

\begin{pf}
Notice that
\[
\int x^2 \,d\mu_{C} = \frac{ 1} { n} \TR\bigl(
C^2 \bigr) \leq\frac{2} { n^2} \sum_{ 1 \leq i \leq j \leq n }
|X_{ij}|^2 \IND\bigl( \veps(n) n^{1/2} <
|X_{i j
} |\leq\veps(n)^{-1} n^{1/2}\bigr).
\]
We may repeat the argument in the proof of Lemma~\ref{le:moder}. This
time we take
$\l= \frac{1}2c \veps(n)^{2 - \alpha}n^{1 + \alpha/ 2}$, where $c$ is
as in \eqref{conx}, and then define
$f$ as in \eqref{eq:deff}. For any $s>0$, one has
\[
\dP{{(\mu_{C} \notin K_{2s} )}} \leq e^{-\l s }
\biggl(1 + e^{ f(
\veps
(n) \sqrt n ) } + \frac{1}2c n^{\alpha/2}
\veps(n)^{-\a}\hspace*{-1pt} \max_{ x
\in
[\veps(n) n^{1/2}, \veps(n)^{-1} n^{1/2}]} e^{f(x)}
\biggr)^{n^2}\hspace*{-1pt}.
\]
Simple considerations show that
$f(x)$, for $x\in[ \veps(n) n^{1/2}, \veps(n)^{-1} n^{1/2}]$ is
maximized at $x = \veps(n) n^{1/2}$, where it satisfies
$f(\veps(n) n^{1/2})\leq-\frac{1}2c\veps(n)^\a n^{\a/2}$. This gives,
for $n$ large enough,
\[
\frac{1}{n^{1 + \alpha/ 2}} \log\dP{{(\mu_{C} \notin K_{2s} )}}
\leq - \frac{1}2 c s \veps(n)^{2 - \alpha} + n^{1- \alpha/ 2}
e^ {-({1}/4)
c \veps(n)^{\alpha} n ^{\alpha/ 2} }.
\]
We choose finally $s = 1 / ( 2 \veps(n)^2 )$. For our choice of $\veps
(n)$ in \eqref{eq:seqveps}, this implies the first claim.

For the second claim, we have
\[
\dP{{\bigl(|I | \geq2 \delta n^{1 + \alpha/2} \bigr)}} \leq\dP \biggl( \sum
_{
1 \leq i \leq j \leq n } \IND \bigl( |X_{ij}| \geq(\log
n)^{2/\alpha} \bigr) \geq\delta n^{1 + \alpha/2} \biggr).
\]
The Bernoulli variables $\IND( |X_{ij}| \geq(\log n)^{2/\alpha} ), 1
\leq i \leq j \leq n$, are independent. Also, by Assumption~\ref{A1},
their average $p_{ij} =\dP(|X_{ij}| \geq(\log n)^{2/\alpha})$ satisfies
\[
p_{ij} \leq p'(n): = e^{ - c (\log n)^2 }
\]
for some $ c > 0$. We argue as in the proof of Lemma~\ref{le:very}.
From Bennett's inequality~\eqref{Bennett},
\begin{eqnarray*}
&&\dP \biggl(\sum_{ 1 \leq i \leq j \leq n } \bigl(\IND \bigl(
|X_{ij}| \geq (\log n)^{2/\alpha} \bigr) - p_{ij}
\bigr) \geq\delta n^{1 + \alpha/2} \biggr) \\
&&\qquad\leq\exp \biggl( - c_0
n^{1 + \alpha/2} \log \biggl( \frac{n^{ \alpha
/2 -1} } { p'(n) } \biggr) \biggr)
\end{eqnarray*}
for some constant $c_0 =c_0(\d)> 0$. Since $p'(n) = o (n^{\alpha
/2-1})$, this gives the claim.
\end{pf}

\subsection{Auxiliary estimates}
To complete the proof of Proposition~\ref{prop:exptight}, we shall need
two extra results. The first is due to Guionnet and Zeitouni \cite{MR1781846},
Corollary~1.4.

\begin{theorem}[(Concentration for matrices with bounded entries)]\label{th:GZ}
Let $\kappa\geq1$, let $Y \in\cH_n(\dC)$ be a random matrix with
independent entries\break $(Y_{ij})_{1 \leq i \leq j \leq n}$ bounded by
$\kappa$, and let $M \in\cH_n(\dC)$ be a deterministic matrix such
that\vspace*{1pt} $\int x^2 \,d\mu_{M} \leq\kappa^2$.
There exists a universal constant $c>0$ such that for all $( c \kappa^2
/ n )^{2/5} \leq t \leq1$,
\[
\dP{{\bigl(W_1 {{(\mu_{ Y / \sqrt n + M },\dE\mu_{ Y / \sqrt n + M } )}}
\geq t \bigr)}} \leq\frac{c \kappa} {t^{ 3/2}} \exp{{\biggl(- \frac{n^2 t^5}{ c \kappa^4 }
\biggr)}}.
\]
\end{theorem}
In \cite{MR1781846}, Corollary~1.4, the result is stated for matrices
$Y$ in $\cH_n (\dC)$ such that the entries have independent real and
imaginary parts. The extension to our setting follows by using a
version of Talagrand's concentration inequality for independent bounded
variables in $\dC$. Also, the matrix $M$ is not present in \cite
{MR1781846}. It is, however, not hard to check that its presence does
not change the argument in \cite{MR1781846}, page 132, since one can
use the bound
\[
\int x^2 \,d\mu_{ Y / \sqrt n + M} \leq2 \int x^2 \,d
\mu_{ Y / \sqrt n} + 2 \int x^2 \,d\mu_{ M } \leq4
\kappa^2.
\]
The latter is an easy consequence of, for example, Lemma~\ref{le:HW}.

The second result we need is
a uniform bound on the rate of the convergence of the empirical
spectral measure of sums of random matrices.

\begin{theorem}[(Uniform asymptotic freeness)]\label{th:UAF}
Let $Y = (Y_{ij} )_{1 \leq i, j \leq n} \in\cH_n (\dC)$ be a Wigner
random matrix with $\VAR( Y_{12} ) = 1$, $\dE|Y_{12} |^3 < \infty$
and $\dE|Y_{11} |^2< \infty$. There exists a universal constant $c >
0$ such that for any integer $n \geq1$ and any $M \in\cH_n (\dC)$,
\[
d {{(\dE\mu_{Y / \sqrt n + M}, \mu_{\mathrm{sc}} \boxplus\mu_M
)}} \leq c \frac{ \sqrt{ \dE|Y_{11} |^2 } + \dE|Y_{12} |^3 }{ \sqrt n }.
\]
\end{theorem}
%
A striking point of the above theorem is that the constant $c$ does not
depend on~$M$. The result is a variation around Pastur and Shcherbina
\cite{MR2808038}, Theorem~18.3.1. The detailed proof of Theorem~\ref
{th:UAF} is given in Appendix~\ref{sec:UAF} below. We are now ready to
finish the proof of Proposition~\ref{prop:exptight}.

\subsection{Proof of Proposition \texorpdfstring{\protect\ref{prop:exptight}}{2.1}}
By Lemmas \ref{le:very} and \ref{le:moder}, it is sufficient to prove that
$
\mu_{A + C}
$ and $\mu_{\mathrm{sc}} \boxplus\mu_{C}$ are exponentially equivalent: for any
$\delta> 0$,
%
\begin{equation}
\label{expeq1} \lim_{n \to\infty} \frac{1}{n^{1 + \alpha/ 2}} \log\dP{{\bigl(d(
\mu _{\mathrm{sc}} \boxplus\mu_{C}, \mu_{A + C} ) \geq
\delta\bigr)}} = - \infty.
\end{equation}
%
Let $\cF$ be the $\si$-algebra generated by the random variables
\[
\{ X_{ij}\IND_{|X_{ij}|\geq(\log n)^{2/\alpha}} \}.
\]
Then the random matrix $C$ is $\cF$-measurable.
Define the event
\[
E = {{\biggl\{\int x^2 \,d\mu_{C} \leq(\log n
)^2 \biggr\}}}.
\]
Then $E \in\cF$. Lemma~\ref{le:exptight} implies that for some
sequence $s_1(n) \to\infty$ and all $n$ large enough,
%
\begin{equation}
\label{eq:expEc} \dP{{\bigl(E^{c}\bigr)}} \leq e^{ -s_1(n) n^{1 + \alpha/ 2}}.
\end{equation}
%
Conditional on $\cF$, $\sqrt n A$ is a random matrix with independent
entries\break $(\sqrt n A_{ij})_{1 \leq i \leq j \leq n}$ bounded by $(\log n
)^{2/\alpha}$. Thus, we may apply Theorem~\ref{th:GZ} with $Y/\sqrt n$
replaced by~$A$, and $M$ replaced by $C$. Using \eqref{eq:boundd} to
replace $W_1(\cdot,\cdot)$ by $d(\cdot,\cdot)$, taking
$t=\d$, and $\kappa= ( \log n )^{2/\alpha}$ in Theorem~\ref{th:GZ},
one has that for all $\d>0$, there is a sequence $s_2(n) \to\infty$,
$n\to\infty$, such that
%
\begin{equation}
\label{eq:expGZ} \IND_E \dP_\cF{{\bigl(d (
\dE_\cF\mu_{A + C}, \mu_{A + C} ) \geq \delta\bigr)}}
\leq e^{ -s_2(n) n^{1 + \alpha/ 2}},
\end{equation}
where $\dP_\cF$ and $\dE_\cF$ are the conditional probability and
expectation given $\cF$.
Notice that Theorem~\ref{th:GZ} can be applied here since on the event
$E$ one has $\int x^2 \,d\mu_{C}\leq(\log n)^2\leq\kappa^2$. Moreover,
\eqref{eq:expGZ} holds uniformly within $E$, since the bound of Theorem~\ref{th:GZ} is uniform with respect to $M$ satisfying $\int x^2 \,d\mu
_{M}\leq\kappa^2$.

From \eqref{eq:expEc} and \eqref{eq:expGZ}, using
the triangle inequality one has that \eqref{expeq1} follows once we
prove that for any $\d>0$:
%
\begin{equation}
\label{expeq2} \lim_{n \to\infty} \frac{1}{n^{1 + \alpha/ 2}} \log\dP{{\bigl(d(
\mu _{\mathrm{sc}} \boxplus\mu_{C}, \dE_\cF
\mu_{A + C} ) \geq\delta\bigr)}} = - \infty.
\end{equation}

Next, we use a coupling argument to remove the dependency between $A$
and $C$. Let $P_n$ be the law of $X_{12}$ conditioned on $\{ |X_{12} |
< (\log n )^{2/\alpha}\}$, and $Q_n$ be the law of $X_{11}$ conditioned
on $\{ |X_{11} | < (\log n )^{2/\alpha}\}$. We also define $I = \{
(i,j) \dvtx |X_{ij} | \geq(\log n ) ^{2/\alpha} \}$. Given $\cF$, if
$(i,j) \in I$, then $A_{ij} = 0$ while, if $(i,j) \notin I$ and $1 \leq
i \leq j \leq n$, then $\sqrt n A_{ij}$ has conditional law $P_n$ or
$Q_n$ depending on whether $i < j$ or $i = j$.

On our probability space, we now consider $Y$ an independent Hermitian
random matrix such that $(Y_{ij})_{1 \leq i \leq j \leq n}$ are
independent, and for $1 \leq i \leq n$, $Y_{ii}$ has law $Q_n$, while
for $1 \leq i < j \leq n$, $Y_{ij}$ has law $P_n$. We form the matrix
\[
A'_{ij} = \IND\bigl( (i, j) \notin I \bigr)
A_{ij} + \IND\bigl( (i, j) \in I \bigr) \frac{ Y_{ij} }{\sqrt n}.
\]
By construction, $\sqrt n A'$ and $Y$ have the same distribution and
are independent of~$\cF$. Also, by Lemma~\ref{le:HW} and Jensen's inequality,
\begin{eqnarray*}
\dE_{\cF} d( \mu_{A + C}, \mu_{A' + C} ) & \leq&\sqrt{
\dE_{\cF
} \frac
{ \TR( A - A')^2} {n } }
\\
& \leq&\sqrt{ \frac{ 1} { n^2} \sum_{ 1 \leq i, j \leq n } \dE
_{\cF} \IND\bigl( (i, j) \in I \bigr) |Y_{ij}|^2}
\leq c_0 \sqrt{ \frac{ | I |} { n^2} },
\end{eqnarray*}
where we have used the fact that,
for some constant $c_0 > 0$,
\[
\max\bigl( \dE|Y_{11} |^2, \dE|Y_{12}|
^2 \bigr) \leq c_0^2.
\]
Define the event
\[
F = {{\bigl\{| I | \leq\delta^2n^2/c_0^2
\bigr\}}}.
\]
Then $F \in\cF$ and
%
\begin{equation}
\label{expeq4} \IND_F\dE_{\cF} d( \mu_{A + C},
\mu_{A' + C} ) \leq\delta.
\end{equation}
From Lemma~\ref{le:exptight}, for some sequence $s_3(n) \to\infty$,
for all $n$ large enough,
%
\begin{equation}
\label{eq:expFc} \dP{{\bigl(F^{c}\bigr)}} \leq e^{ -s_3(n) n^{1 + \alpha/ 2}}.
\end{equation}
Observe that by definition of the distance \eqref{eq:defdist},
\[
d ( \dE_\cF\mu_{A' + C}, \dE_\cF
\mu_{A + C} )\leq\dE_\cF d ( \mu _{A' + C},
\mu_{A + C} ).
\]
Since $A'$ and $Y/\sqrt n $ have the same distribution,
we deduce from \eqref{expeq4}, \eqref{eq:expFc} and the triangle
inequality that
the proof of \eqref{expeq2} can be reduced to the proof of
%
\begin{equation}
\label{expeq5} \lim_{n \to\infty} \frac{1}{n^{1 + \alpha/ 2}} \log\dP{{\bigl(d(
\mu _{\mathrm{sc}} \boxplus\mu_{C}, \dE_\cF
\mu_{ Y/\sqrt n + C} ) \geq\delta \bigr)}} = - \infty.
\end{equation}
%
Clearly, $\dE|Y_{12}| ^3 \leq c_0 (\log n )^{6/\alpha}$ and $\sigma^ 2
= \VAR( Y_{12} ) \to1$. We may apply the uniform estimate of Theorem~\ref{th:UAF}, applied to $Y/(\si\sqrt n)$ and $M = C$, which is $\cF
$-measurable. We find
for any $\d>0$,
\[
\dP{{\bigl(d {{(\mu_{\mathrm{sc}} \boxplus\mu_{C},
\dE_\cF\mu_{ {Y}/ {(\si
\sqrt n)} + C} )}} \geq\delta\bigr)}}=0
\]
for all $n\geq n_0(\d)$ where $n_0(\d)$ is a constant depending only on
$\d$.

On the other hand, arguing as above, from Hoffman--Wielandt's
inequality (Lemma~\ref{le:HW}) and Jensen's inequality, for any
$\delta>0$,
\begin{eqnarray*}
d {{(\dE_{\cF} \mu_{{Y} /{\sqrt n} + C}, \dE_{\cF} \mu
_{{Y}/{ (\si\sqrt n)} + C} )}} & \leq& \dE_{\cF} d
 {{(\mu_{{Y}/ {\sqrt n} + C},
\mu_{{Y}/ {( \si\sqrt n)} + C} )}} \\
&\leq&\dE_\cF \sqrt{ \frac{(1
- {1}/ {\si})^ 2 }{n^2} \TR
\bigl( Y^2 \bigr)}
\\
& \leq& \frac{|1 - {1} /{\si}|}{ n} \sqrt{ \dE\TR\bigl( Y^2 \bigr) } \leq
\delta
\end{eqnarray*}
for all $n \geq n_1(\d)$ where $n_1(\d)$ is a constant depending only
on $\d$.

This concludes the proof of \eqref{expeq5} and of Proposition~\ref
{prop:exptight}.


\begin{remark}\label{weaker}
In the proof of Proposition~\ref{prop:exptight}, we have only used the
following assumptions on the Wigner matrix $X$: (i) $\VAR(X_{12}) =1$
and (ii) there exists $c >0$ such that for all $i\leq j$,
\[
\limsup_{t\to\infty}\frac{1}{t^{\alpha}} \log\dP{{\bigl(|X_{ij}|
\geq t \bigr)}}\leq-c. 
\]
\end{remark}

\section{Large deviations of very sparse rooted networks}
\label{sec:esgraphs}
In this section, we start by adapting to our setting the notion of
local weak convergence of rooted networks, introduced in \cite{bensch,aldoussteele} and \cite{aldouslyons}. Next, we introduce a
suitable projective limit topology on the space of networks. Then we
prove the LDP for the network $G_n$ induced by the very sparse matrix
$C$. Finally, we introduce the spectral measure associated to a network
and project the LDP for networks onto a LDP for spectral measures.

\subsection{Locally finite Hermitian networks}
Let $V$ be a countable set, the \emph{vertex} set. A pair $(u,v)\in V^2$
is an \emph{oriented edge}.
A \emph{network} or \emph{weighted graph} $G = (V, \omega)$ is a vertex
set $V$ together with a map
$\omega$ from $V^2$ to $\dC$.
We say that a network is \emph{Hermitian}, if for all $(u,v) \in V^2$,
\[
\omega(u,v) = \overline{ \omega(v,u)}.
\]
For ease of notation, we sometimes set $\omega(v) = \omega(v,v)$ for
the weight of the loop at $v$.
The degree of $v$ in $G$ is defined by
\[
\deg(v) = \sum_{u \in V } \bigl|\omega(v,u)\bigr |^2.
\]
The network $G$ is \emph{locally finite} if for any vertex $v$, $\deg
(v)< \infty$.

A path $\pi$ from $u$ to $v$ in $V$ is a sequence $\pi= (u_0, \ldots,
u_k)$ with $u_0 = u$, $u_k = v$ and, for $1 \leq i \leq k$, $| \omega(
u_{i-1}, u_{i} ) | > 0$. If such $\pi\dvtx u\to v$ exists, then one defines
the $\ell_2$ distance
\[
D_\pi(u,v)= \Biggl( \sum_{i=1}
^k\bigl| \omega( u_{i-1}, u_{i} ) \bigr|^{-2}
\Biggr) ^{1/2}.
\]
The distance between $u$ and $v$ is defined as
\[
D(u,v)=\inf_{\pi: u\to v} D_\pi(u,v).
\]
Notice that weights are thought of as inverse of distances. If there is
no path $\pi\dvtx u\to v$, then the distance $D(u,v)$ is set to be infinite.
A network is \emph{connected} if $D(u,v)<\infty$ for any $u\neq v \in V$.

All networks we consider below will be Hermitian and locally finite,
but not necessarily connected. We call $\cG$ the set of all such networks.
For a network $G\in\cG$, to avoid possible confusion, we will often
denote by $V_G$, $\omega_G$, $\deg_G$ the corresponding vertex set,
weight and degree functions.

Clearly, any $n\times n$ Hermitian matrix $H_n\in\cH_n(\dC)$ defines a
finite network $G=G(H_n)$ in a natural way, by taking
%
\begin{equation}
\label{ghn} V_{G}=\{1,\ldots,n\},\qquad \o_{G}(i,j)=H_n(i,j)
.
\end{equation}
For simplicity, we often write simply $H_n$ instead of $G(H_n)$.

\subsection{Rooted networks}
Below, a \emph{rooted network} $(G,o) = ( V, \omega, o)$ is a Hermitian,
locally finite and connected network $(V, \omega)$
with a distinguished vertex $o \in V$, the root. For $t > 0$, we denote
by $(G,o)_t$ the rooted network with vertex set $\{u\in V\dvtx D(o,u)\leq
t\}$, and with the weights induced by $\o$.
Two rooted networks $(G_i,o_i) = ( V_i, \omega_i, o_i )$,
$i \in\{1,2\}$, are \emph{isomorphic} if there exists a bijection
$\si
\dvtx V_{1} \to V_{2}$ such that $\si( o_1) = o_2$ and $\si( G_1) = G_2$,
where $\si$ acts on $G_1$ through $\si( u, v ) = (\si( u), \si(v)
)$ and $\si( \omega) = \omega\circ\si$.

We define the semidistance $d_\loc$ between two rooted networks
$(G_1,o_1)$ and $(G_2,o_2)$ to be
\[
d_{\loc} \bigl((G_1,o_1),( G_2,o_2)
\bigr) = \frac{1}{1 + T},
\]
where $T$ is the supremum of those $t > 0$ such that
there is a bijection $\si\dvtx   V_{(G_1,o_1)_t} \to V_{(G_2,o_2)_t}$ with
$\si(o_1) = o_2$ and such that the function $ \omega_{G_2} \circ\si-
\omega_{G_1}$ is bounded by $1/t$ on $V^2_{(G_1,o_1)_t}$.

The rooted network isomorphism defines a space $\cG_*$ of equivalence
classes of rooted networks $(G,o)$.
On the space $\cG_*$, $\dloc$ becomes a distance. The associated
topology will be referred to as the \emph{local topology}.
We write $\frg$ for an element of $\cG_*$.
We shall denote the convergence on $(\cG_*,\dloc)$ by $\dloc(\frg
_n,\frg
)\to0$ or $\frg_n \locto\frg$.

The space $(\cG_*,\dloc)$ is separable and complete \cite{aldouslyons}.
Let $\cP(\cG_*)$ denote the space of probability measures on $\cG_*$.
For $\mu,\mu_n\in\cP(\cG_*)$, we write $\mu_n \locweak\mu$ when
$\mu
_n$ converges weakly, that is, when $\int f \,d\mu_n \to\int fd\mu$ for
every bounded continuous function $f$ on $(\cG_*,\dloc)$. This notion
of weak convergence is often referred to as \emph{local weak
convergence}. See \cite{aldouslyons} for more details and examples.

For a network $G\in\cG$, and $v\in V_G$, one writes $G(v)$ for the
connected component of $G$ at $v$, that is, the largest connected
network $G'\subset G$ with $v\in V_{G'}$.
If $G\in\cG$ is finite, that is, $V_G$ is finite,
one defines the probability measure $U(G)\in\cP(\cG_*)$ as the law of
the equivalence class of the rooted network
$(G(o),o)$ where the root $o$ is sampled uniformly at random from $V_G$:
\[
U(G)=\frac{1}{V_G}\sum_{v\in V_G}
\d_{\frg(v)},
\]
where $\frg(v)$ stands for the equivalence class of $(G(v),v)$.
If $G_n, n \geq1$, is a sequence of finite networks from $\cG$, we
shall say that $G_n$ has \emph{local weak limit} $\rho\in\cP(\cG
_*)$ if
$U(G_n) \locweak\rho$.

\subsection{Sofic measures} Following \cite{aldouslyons}, a measure
$\r
\in\cP(\cG_*)$ is called \emph{sofic} if there exists a sequence of
finite networks $G_n, n \geq1$, whose local weak limit is~$\rho$.
We shall need to identify a subset of the sofic measures.
Let $\vartheta_a,\vartheta_b$
denote the laws of $X_{12}/|X_{12}|$
and $X_{11}/|X_{11}|$,
respectively, for
$X_{12}\in\cS_\a(a)$ and $X_{11}\in\cS_\a(b)$; see Assumption~\ref
{A1}, and let $S_a, S_b\subset\dS^1$ denote their supports.
Let $\cA_n\subset\cH_n(\dC)$ be the set of $n\times n$ Hermitian
matrices $H$ such that
either $H_{ij}=0$ or $H_{ij}/|H_{ij}|\in S_a$ for all $i<j$, and such
that either $H_{ii}=0$ or $H_{ii}/|H_{ii}|\in S_b$ for all $i$.
We say that $\r\in\cP(\cG_*)$ is \emph{admissible sofic} if there exists
a sequence of matrices
$H_n\in\cA_n$ such that $U(H_n)\locweak\r$, where $H_n$ is identified
with the associated network $G(H_n)$ as in \eqref{ghn}.
We denote by $\cP_s ( \cG_*)$ the set of admissible sofic probability
measures.
Measures in $\cP_s ( \cG_*)$ will often be called simply sofic if no
confusion can arise.

Let $\frg_\varnothing$ stand for the trivial network consisting of a
single isolated vertex (the root)
with zero weights. We refer to $\frg_\varnothing$ as the empty network.
Clearly, the Dirac mass at the empty network $\r=\d_{\frg_\varnothing}$
is sofic
(it suffices to consider matrices with zero entries). Let us consider
some more examples.

\begin{example}\label{sinet}
Suppose that $S_b=\{-1,+1\}$. Let
$Y_1,Y_2,\ldots$ be i.i.d. random variables with distribution $\nu\in
\cP
(\dR)$.
Consider the random diagonal matrix $H_n$ with $H_n(i,i)=Y_i$. Then, by
the law of large numbers,
almost surely
$U(H_n)\locweak\r$, where $\r$ is given by
\[
\r=\int_{\dR}\d_{{\frg}_{x}}\,d\nu(x),
\]
if ${\frg}_{x}$ is the network consisting of a single vertex (the root)
with loop weight equal to $x$.
\end{example}

\begin{example}\label{dinet}
Suppose that $Z_1,Z_3,Z_5,\ldots$
are i.i.d. complex random variables with law $\mu\in\cP(\dC)$ such
that $\mu$-a.s. one has either $Z_1=0$, or $Z_1/|Z_1|\in S_a$.
Consider the $n\times n$ matrix
$H$
such that $H_n(j,j+1)=Z_j$, $H_n (j+1,j)=\bar{Z}_j$,
for all odd $1\leq j\leq n-1$, and all other entries of $H_n$ are zero.
By construction, $H_n\in\cA_n$ almost surely.
From the law of large numbers, almost surely $U(H_n)\locweak\r$,
where $\r$ is given by
\[
\r=\frac{1}2\int_{\dC}(\d_{\hat{\frg}_z}+
\d_{\hat{\frg}_{\bar
z}})\,d\mu(z),
\]
if $\hat{\frg}_{z}$ denotes the equivalence class of the two vertex
network $(V,\o,o)$,
with $V=\{o,1\}$, $\o(o,1)=z$, $\o(1,o)=\bar z$ and $\o(o,o)=\o(1,1)=0$.
\end{example}

\begin{example}\label{nnet}
For any fixed $n\in\dN$, if $H_n\in\cA_n$, then $U(H_n)\in\cP
_s (
\cG_*)$.
Indeed, take a sequence of $m\times m$ matrices $A_m\in\cA_m$ defined
as follows. Let $k,r\geq0$, with $r<n$, be integers such that
$m=kn+r$, and take $A_m$ as the block diagonal matrix with the first
$k$ blocks all equal to
$H_n$ and the last block of size $r$ equal to zero. Then $U(A_m)=\frac
{n}{n+(r/k)} U(H_n) + \frac{1}{1+(kn/r)}\d_{\frg_\varnothing}$. As
$m\to
\infty$, $r/k\to0$, $kn/r\to\infty$ and, therefore, $U(A_m)$ converges
to $U(H_n)$.
\end{example}

\subsection{Truncated networks}
It will be important to work with suitable truncations of the weights.
To this end we consider,
for $0 < \theta< 1$, networks $G\in\cG$ such that
for any $(u,v) \in V_G^2$,
%
\begin{equation}
\label{thetanet} \deg_G ( v) \leq\theta^{-2} \quad\mbox{and}\quad
\bigl| \omega_G (u,v) \bigr| \geq\theta\IND\bigl(\omega_G (u,v)
\neq0\bigr).
\end{equation}
We call $\cG^\theta$ the set of all such networks. Clearly, any $G\in
\cG
^\theta$ is locally finite and has at most $\theta^{-4}$ outgoing
nonzero edges from any vertex.
As before, one defines the space $ \cG^{\theta}_*$ by taking
equivalence classes of connected rooted networks from $\cG^\theta$.
We define $\cP( \cG^{\theta}_*)$ as the sets of $\rho\in\cP( \cG_*)$
with support in $\cG^\theta_*$, and set $\cP_s( \cG^{\theta
}_*)=\cP( \cG
^{\theta}_*)\cap\cP_s( \cG_*)$.

\begin{lemma}
\label{le:LWprop}
\textup{(i)} 
$\cP_s (\cG_*)$ is closed for the local weak topology.
\begin{longlist}[(ii)]
\item[(ii)]
For any $\theta> 0$, $\cG^{\theta}_*$ is a compact set for the local
topology.
\end{longlist}
\end{lemma}

\begin{pf}
For (i): by definition, $\cP_s(\cG^*)$ is the closure of
the set of $U(G)$ such that $G$ is an admissible finite network [i.e.,
for some integer $n\geq1$, $H\in\cA_n$ and $G = G(H)$ as in \eqref{ghn}].

For (ii): let $\frg\in\cG^{\theta}_*$ and $(G,o)$ be a
rooted network in the equivalence class $\frg$. Observe that each edge
of $G$ has a weight bounded above by $\theta^{-1}$. This implies that
in $G$ each path whose total length is bounded by $t>0$, contains at
most $t^2 / \theta^2$ edges. Moreover, $G$ has at most $\theta^{-4}$
outgoing edges from any vertex. Hence, $G$ has at most $n (t) = \theta
^{-4 t^2 / \theta^2}$ vertices at distance less than $t$ from any
given vertex.

Now, we denote by $\cG^{\theta,t}_*$ the set of equivalence classes of
$(G,o)_t$ such that the equivalence class of $(G,o)$ is in $\cG
^{\theta
}_*$. There is a finite number, say $m(t)$, of equivalence classes of
rooted connected graphs with less than $n(t)$\vadjust{\goodbreak} vertices (without
weights). Since all weights of $\frg\in\cG^\theta_*$ are in
$[\theta
,\theta^ {-1}]$, there is a covering of $\cG^{\theta,t}_*$ with balls
of radius $1/(1+t)$ of cardinal at most $k(t) = m(t) ( t \theta^ {-1}
)^ {n(t)^ 2}$.

Notice that for any rooted network $d_{\loc} ( (G,o), (G,o)_t ) \leq
1/(1+t)$. Hence, from the definition of $d_{\loc}$, we have proved
that, for any $t>0$, there exists a finite covering of $ \cG^{\theta
}_*$ with balls of radius $1/(1+t)$. This proves that $ \cG^{\theta}_*$
is precompact. The fact that $\cG^{\theta}_*$ is closed follows
directly from \eqref{thetanet}.
\end{pf}

Next, we describe a canonical way to obtain a network in $\cG^\theta$
by truncating a network from $\cG$.
This will allow us to introduce a topology on $\cP(\cG_*)$ that is
weaker than the local weak topology.
In particular, a topology for which $\cP_s (\cG_*)$ is compact; compare
Lemmas \ref{le:LWprop} and~\ref{le:LWprop0}.
For $0 < \theta< 1$, define the two continuous functions
\begin{eqnarray*}
\chi_\theta(x) &=& \cases{ %
0, & \quad$\mbox{if }  x \in[0,\theta),$
\vspace*{2pt}\cr
( x - \theta) / \theta,& \quad$\mbox{if }  x \in[\theta,2 \theta),$
\vspace*{2pt}\cr
1, & \quad$\mbox{if }  x \in[2\theta, \infty),$}
\\
 \wchi_\theta(x) &=& \cases{
 1, &\quad
$\mbox{if }  x \in\bigl[0,\theta^{-2} -1 \bigr),$
\vspace*{2pt}\cr
\theta^{-2} - x, & \quad$\mbox{if }  x \in\bigl[\theta^{-2} -1,
\theta^{-2}\bigr),$
\vspace*{2pt}\cr
0, & \quad$\mbox{if }  x \in\bigl[\theta^{-2}, \infty\bigr)$}
\end{eqnarray*}
that will serve as
approximations for the indicator functions $\IND( x \geq\theta)$ and
$\IND( x \leq\theta^{-2})$.

If $G = (V,\omega) $, we define $\wt G_\theta= (V,\wt\omega_\theta)$
as the network with vertex set $V$ and, for all $u,v \in V$,
%
\begin{equation}
\label{eq:defomthet} \wt\omega_\theta( u,v) = \omega( u, v )
\wchi_\theta{{\bigl(\deg_{G} (u) \vee\deg_{G} (v)
\bigr)}}.
\end{equation}
Next, we define $G_\theta= (V,\omega_\theta)$ as the network with
vertex set $V$ and, for all \mbox{$u,v \in V$},
%
\begin{equation}
\label{eq:defomthet2} \omega_\theta( u,v) = \wt\omega_\theta( u, v )
\chi_\theta {{\bigl(\bigl|\wt\omega_\theta( u, v )\bigr| \bigr)}}.
\end{equation}
Clearly, $G_\theta$ satisfies \eqref{thetanet}, and for any $u, v \in
V$, 
$|\omega_{G_\theta} ( u, v) |\leq\theta^{-1}$, and
%
\begin{equation}
\label{monot} \deg_{G_\theta} ( u) \leq\deg_{G} ( u) \quad\mbox{and}
\quad \bigl|\omega_{G_\theta} ( u, v)\bigr |\leq\bigl|\omega_{G} ( u, v)\bigr |.
\end{equation}
If $\frg\in\cG_*$ and the network $(G,o)$ is in the equivalence class
$\frg$, then $\frg_\theta\in\cG_*^\theta$ is defined as the
equivalence class of $(G_\theta(o), o)$, where $G_\theta$ is defined
by \eqref{eq:defomthet2}. This defines a map $\frg\mapsto\frg
_\theta$
from $\cG_*$ to $\cG_*^\theta$. If $\rho\in\cP( \cG_*)$ and
$\frg$
has law $\rho$, the law of $\frg_\theta$ defines a new measure $\rho
_\theta\in\cP( \cG^\theta_*)$.

The next lemma follows easily from the continuity of $\chi_\theta
,\wchi
_\theta$ and the fact that as $\theta\to0$, for any for $x > 0$,
$\chi
_\theta(x) \to1$ and $\wchi_\theta(x) \to1$.

\begin{lemma}[(Continuity of projections)]\label{le:contproj}
\begin{longlist}[(iii)]
\item[(i)] For $\theta> 0$, the map $\frg\mapsto\frg_\theta$ from $\cG_*
\to\cG_*^\theta$ is continuous for the local topology;\vadjust{\goodbreak}
\item[(ii)]
for $\theta> 0$, the map $\rho\mapsto\rho_\theta$ from $\cP(\cG_*)$
to $\cP(\cG_*^\theta)$ is continuous for the local weak topology;
\item[(iii)]
as $\theta\to0$, one has $\frg_\theta\locto\frg$ and $\rho
_\theta
\locweak\rho$, for any $\frg\in\cG_*$ and $\rho\in\cP(\cG_*)$.
\end{longlist}
\end{lemma}

\subsection{Projective topology for locally finite rooted networks}
In order to circumvent the lack of compacity of $\cP_s ( \cG_*)$
w.r.t. local weak topology, we now introduce a weaker topology, the
\emph{projective topology}.
For integers $j \geq1$, set
\[
\theta_j = 2^{-j}.
\]
Let $p_{j}\dvtx \cG_*\to\cG_*^{\theta_j}$ be defined by $p_j ( \frg) =
\frg
_{\theta_j}$.
Similarly, for $1 \leq i \leq j$, $p_{ij}\dvtx \cG_* ^{\theta_j}\to\cG_*
^{\theta_i}$ is the map $p_{ij}(\frg)=\frg_{\theta_i}$, $\frg\in
\cG_*
^{\theta_j}$.
The collection $(p_{ij})_{ 1 \leq i \leq j}$ is a \emph{projective
system} in the sense that for any $1 \leq i \leq j \leq k$,
%
\begin{equation}
\label{eq:const} p_{i k } = p_{ij} \circ p_{jk}.
\end{equation}
The latter follows from $2 \theta_{j+1} \leq\theta_{j}$ and $\theta
^{-2}_j \leq\theta^{-2}_{j+1} - 1$.

Define the projective space $\wcG_* \subset\prod_{ j \geq1} \cG
^{\theta_j}_* $ as the set of $y=(y_1,y_2,\ldots) \in\prod_{ j \geq1}
\cG^{\theta_j}_*$ such that for any $i \leq j$, $p_{ ij} ( y_j) = y_i$;
see, for example, \cite{dembo}, Appendix~B,
for more details on projective spaces. One can identify $\cG_*$ and~$\wcG_*$:

\begin{lemma}
\label{le:defproj}
The
map $\iota( \frg) =(p_j ( \frg))_{j \geq1}$ from $\cG_*$ to $\wcG_*$
is bijective.
\end{lemma}

\begin{pf}
The fact that $\iota$ is injective is a consequence of Lemma~\ref
{le:contproj} part (iii). It remains to prove that the map $\iota$ is
surjective. Let $y = (y_j) \in\wcG_*$.
One can represent the $y_j$'s by rooted networks $(G_j,o)= (V_j, \omega
_j, o )$ such that $V_j\subset V_{j+1}$. Set $V:=\bigcup_{j\geq1}V_j$.
By adding isolated points, one can view $(G_j,o)$ as the connected
component at the root of the network $\hat G_j=(V,\o_j)$, where $\o
_j(u,v)=0$ whenever either $u$ or $v$ (or both) belong to $V\setminus
V_j$. Moreover, one has that $\hat G_i=(\hat G_j)_{\theta_i}$ for all
$i<j$. This sequence of networks is monotone in the sense of \eqref{monot}.

For fixed $u, v \in V$, and $j\in\dN$,
if $\o_j(u,v)\neq0$ then the degree of $u$ and $v$ is bounded by
$2^{2j}$ in any network $\hat G_k$, $k\geq j$ and, therefore, $\o
_k(u,v)=\o_{j+1}(u,v)$ for all $k\geq j+1$. In particular, for all $u,
v \in V$
the limit
\[
\omega(u,v) = \lim_{ j \to\infty} \omega_{j} (u,v)
\]
exists and is finite. The same argument shows that
for any $u\in V$,\break $\lim_{ j \to\infty} \deg_{\hat G_j} (u)$ exists and
equals
\[
\sum_{v \in V} \bigl| \omega(u,v)\bigr| ^2 < \infty.
\]
To prove surjectivity of the map\vspace*{1pt} $\iota$, it suffices to take the
network $G=(V,\o)$, and observe that it satisfies $G_{\theta_j}=\hat
G_j$ for all $j\in\dN$.
\end{pf}

With a slight abuse of notation, we will from now on write $\cG_*$ in
place of $\wcG_*$. The projective topology on $\cG_*$ is the topology
induced by the metric
\[
\dproj\bigl( \frg, \frg'\bigr) = \sum
_{ j \geq1} 2^{-j} \dloc\bigl( \frg _{\theta
_j},
\frg'_{\theta_j} \bigr).
\]
The metric space $(\cG_*,\dproj)$ is complete and separable. Also,
$\frg_n \projto\frg$, that is, $\dproj(\frg_n,\frg)\to0$, if and
only if for any $\theta> 0$, $(\frg_n)_{\theta} \locto\frg_\theta$.
The \emph{projective weak topology} is the weak topology on $\cP( \cG
_*) $ associated to continuous functions on $(\cG_*, \dproj)$. We
denote the associated convergence by $\projweak$. Notice that $\rho_n
\projweak\rho$ if and only if for any $\theta> 0$, $(\rho
_n)_{\theta
}\locweak\rho_\theta$.
The topology generated by $\dproj$ is coarser than the topology
generated by $\dloc$, and the weak topology associated to $\projweak$
is coarser than the weak topology associated to $\locweak$.

\begin{example}\label{starex}
Consider the star shaped rooted network $(G_n,1) =\break (V_n, \omega_n
,1)$ where $V_n = \{1, \ldots, n\}$, with $\o_n(u,v)=\o_n(v,u)=1$, if
$u=1$ and $v\neq1$, and $\o(u,v)=0$ otherwise.
Let $\frg_n$ denote the associated equivalence class in $\cG_*$.
Then $\frg_n$ does not converge in $(\cG_*, \dloc)$ because of the
diverging degree at the root. However, in $(\cG_*, \dproj)$, $\frg_n
\projto\frg_\varnothing$ where $\frg_\varnothing$ is the empty network.
Moreover, $U(G_n)$ does not converge in $\cP(\cG_*)$ for $ \locweak
$ however $U(G_n) \projweak\delta_{\frg_\varnothing}$.
\end{example}

\begin{lemma}
\label{le:LWprop0}
\textup{(i)} 
$\cG_*$ is compact for the projective topology.
\begin{longlist}[(ii)]
\item[(ii)]
$\cP_s (\cG_*)$ is compact for the projective weak topology.
\end{longlist}
\end{lemma}

\begin{pf} Statement (i) is a consequence of Tychonoff
theorem and Lem\-ma~\ref{le:LWprop}(ii). It implies that
$\cP( \cG_*)$ is compact for projective weak topology. Hence, to prove
statement (ii), it is sufficient to check that $\cP_s
(\cG_*)$ is closed. Assume that $\rho_n \in\cP_s(\cG_*)$ and $\rho_n
\projweak\rho$. Then for any $\theta> 0$, $(\rho_n)_\theta\in\cP_s
( \cG_*)$ and $(\rho_n)_{\theta} \locweak\rho_\theta$. By Lemma~\ref
{le:LWprop}(i), we deduce that $\rho_\theta\in\cP
_s(\cG_*)$.
However, as $\theta\to0$, using Lemma~\ref{le:contproj}, we find
$\rho
_\theta\locweak\rho$. By appealing to Lemma~\ref{le:LWprop}(i) again, we get $\rho\in\cP_s( \cG_*)$.
\end{pf}

\subsection{Large deviations for the network $G_n$}\label{ldpnetw}
For a rooted network $(G,o)$, $G=(V_G,\o_G)$, define the functions
%
\begin{equation}
\label{psiphi} \psi(G,o) = \bigl| \omega_G ( o)\bigr |^\alpha\quad\mbox{and}\quad \phi (G,o) = \frac{1} 2 \sum_{ v \in V_G \setminus o } \bigl|
\omega_G ( o, v) \bigr|^\alpha.
\end{equation}
Since these functions are invariant under rooted isomorphisms, one can
take them as functions on
$\cG_*$. Then, if $\r\in\cP(\cG_*)$ we write $\dE_\r\psi$, and
$\dE_\r
\phi$ to denote the corresponding expectations. We remark that for any
$\theta>0$, the restriction of $\phi,\psi$ to $(\cG_*^\theta,\dloc)$
gives two bounded continuous functions. Therefore, as functions on
$(\cG
_*,\dproj)$, $\phi$ and $\psi$ are lower semicontinuous.

We now come back to the random matrix $C=C(n)$ defined in \eqref
{decompo}. For integer $n\geq1$, consider the associated network
%
\begin{equation}
\label{network} G_{n}=(V_n,\o_n)\qquad \mbox{with }
V_n = \{ 1, \ldots, n\} \mbox{ and } \omega_n(i,j) =
C_{ij}.
\end{equation}
From the first Borel--Cantelli lemma, almost surely the matrix $C$ has
no nonzero entry for $n$ large enough. Therefore, almost surely,
$U(G_n) \locweak\d_{\frg_\varnothing}$,
the Dirac mass at the empty network.

For ease of notation, we define the random probability measure
\[
\rho_n = U(G_n).
\]
Notice that, by definition one has
%
\begin{equation}
\label{ern1} \dE_{ \rho_n } \psi= \frac{1} {n^{1+\alpha/2} } \sum
_{i=1} ^n |X_{ii}|^{\alpha} \IND
\bigl( \veps(n) \sqrt n \leq| X_{ii} | \leq\veps (n)^{-1}
\sqrt n \bigr)
\end{equation}
and
%
\begin{equation}
\label{ern2} \dE_{ \rho_n } \phi= \frac{1} {n^{1+\alpha/2} } \sum
_{1 \leq i < j
\leq n} |X_{ij}|^{\alpha} \IND\bigl( \veps(n)
\sqrt n \leq| X_{ij} | \leq \veps(n)^{-1} \sqrt n \bigr).
\end{equation}
The next proposition gives the large deviation principle for $\r
_n=U(G_n)$ for the projective weak topology.

\begin{proposition}
\label{prop:LDPnetwork}
$U(G_n)$ satisfies an LDP on $\cP(\cG_*)$ equipped with the projective
weak topology,
with speed $n^{1 + \alpha/2}$ and good rate function $I\dvtx \cP(\cG
_*)\mapsto[0,\infty]$ defined by 
%
\begin{equation}
\label{Irho} I (\rho) = \cases{ b \dE_\r\psi+ a
\dE_\r\phi,& \quad$\mbox{if } \r\in\cP_s (\cG_*),$\vspace*{2pt}
\cr
+\infty,&\quad $\mbox{if } \r\notin\cP_s (\cG_*)$.}
\end{equation}
If $a$ or $b$ is equal to $\infty$, the above formula holds with the
convention $\infty\times0 = 0$.
\end{proposition}

\begin{pf} By construction, $\r_n=U(G_n)\in\cP_s( \cG_*)$; see Example~\ref
{nnet}. 
Since $\cP_s( \cG_*)$ is closed (see Lemma~\ref{le:LWprop}), it is
sufficient to establish the LDP on the space $\cP_s( \cG_*)$ with good
rate function $I(\r)=b \dE_\r\psi+ a \dE_\r\phi$, $\r\in\cP_s
(\cG_*)$.

Let $B_\proj(\rho,\delta)$ [resp., $B_\loc(\rho,\delta)$] denote the
closed ball with radius $\delta> 0$ and center $\rho\in\cP_s ( \cG
_*)$ for the L\'evy metric associated to the projective weak topology
(resp., local weak topology).

\emph{Upper bound.} By Lemma~\ref{le:LWprop0}(ii), $\cP_s(\cG_*)$ is compact. Hence, it is sufficient to
prove (see, e.g., \cite{dembo}) that for any $\r\in\cP_s(\cG_*)$
%
\begin{equation}
\label{eq:LDPuptbp} \limsup_{\delta\to0} \limsup_{n \to\infty}
\frac{ 1}{ n^{1 +
\alpha
/2} } \log\dP{{\bigl(\rho_n \in B_\proj( \rho,
\delta) \bigr)}} \leq- b \dE_\rho\psi- a \dE_\rho\phi.
\end{equation}

Assume first that $\dE_\rho\psi$ and $\dE_\rho\phi$ are finite.
From standard properties of weak convergence, and the fact that $\phi
,\psi$ are lower semicontinuous on $(\cG_*,\dproj)$, it follows that
the maps $\mu\mapsto\dE_{\mu} \psi$ and $\mu\mapsto\dE_{\mu}
\phi
$ are lower semicontinuous on $\cP_s (\cG_*)$ w.r.t. the projective
weak topology. Hence, we have for some continuous function $h(\cdot)$
with $h(0) = 0$,
\begin{eqnarray*}
\dP\bigl( \rho_n \in B_\proj( \rho, \delta) \bigr)& \leq
\dP{{\bigl(\dE _{\rho
_n} \psi\geq\dE_\rho\psi- h(\delta);
\dE_{\rho_n} \phi\geq \dE_\rho\phi- h(\delta) \bigr)}}.
\end{eqnarray*}
%
Since \eqref{ern1} and \eqref{ern2} are independent random variables,
%
\begin{eqnarray}\label{eq:lbLDPdede}
&& \dP\bigl( \r_n \in B_\proj( \rho, \delta) \bigr)
\nonumber
\\[-8pt]
\\[-8pt]
\nonumber
&&\qquad \leq\dP{{\bigl(\dE_{ \rho_n } \psi\geq\dE_\rho\psi- h(\delta)
\bigr)}} \dP{{\bigl(\dE_{\rho_n} \phi\geq\dE_\rho\phi- h(\delta)
\bigr)}}.
\end{eqnarray}
To prove the part of the bound involving $\phi$, one may assume $ \dE
_\rho\phi> 0$. Take $\delta$ small enough, so that $s:= \dE_\rho
\phi- h(\delta) > 0$. From \eqref{ern2}, using Markov's inequality,
for any $a_1>0$,
%
\[
\dP{{(\dE_{ \rho_n } \phi\geq s )}} \leq e^{- a_1 n^{1 + \alpha/2}
s } {{\bigl(\dE\exp{{
\bigl(a_1 |X_{12}|^{\alpha} \IND_{ \veps(n) \sqrt n
\leq| X_{12} | \leq\veps(n)^{-1} \sqrt n }
\bigr)}} \bigr)}}^{n(n-1)/2}.
\]
Take $0< a_1< a$.
By assumption, there exists $ a_2 \in(a_1, a)$, such that for all $t >
0$ large enough,
\[
\dP{{\bigl(|X_{12} | \geq t \bigr)}} \leq\exp\bigl( - a_2
t^\alpha\bigr).
\]
Using \eqref{eq:upIbP}, one deduces that
\begin{eqnarray*}
&&\dE\exp{{\bigl(a_1 |X_{12}|^{\alpha}
\IND_{ \veps(n) \sqrt n \leq|
X_{12} | \leq\veps(n)^{-1} \sqrt n } \bigr)}}
\\
&&\qquad \leq1 + e^{ - (a_2 - a_1) \veps(n)^\alpha n^{\alpha/2}
} + \alpha a_1 \int_{ \veps(n) \sqrt n}
^ { \veps(n)^{-1} \sqrt n } x^{\alpha-1} e^{ - ( a_2 - a_1) x^\alpha} \,dx
\\
& &\qquad\leq1 + \frac{ a_2} { a_2 - a_1 } e^{ - (a_2 - a_1)
\veps(n)^\alpha n^{\alpha/2} }.
\end{eqnarray*}
Therefore,
\[
\dP{{(\dE_{ \rho_n } \phi\geq s )}} \leq\exp{ \biggl(- a_1
n^{1 +
\alpha/2} s+ \frac{ a_2} { 2(a_2 - a_1) } n^2 e^{ - (a_2 - a_1) \veps
(n)^\a n^{\a/2}} \biggr)
}.
\]
We have thus proved that for $\d$ small enough
\[
\limsup_{n \to\infty} \frac{ 1}{ n^{1 + \alpha/2} } \log\dP {{(\dE
_{\rho_n} \phi\geq s )}} \leq- a_1 \bigl(\dE_\rho
\phi- h(\delta) \bigr).
\]
Since the above inequality is true for any $a_1 < a$, it also holds for
$a_1 = a$. Similarly, one has
\[
\limsup_{n \to\infty} \frac{ 1}{ n^{1 + \alpha/2} } \log\dP {{(\dE
_{\rho_n } \psi\geq s )}} \leq- b \bigl(\dE_\rho\psi- h(\delta )
\bigr).
\]
From \eqref{eq:lbLDPdede}, it follows that \eqref{eq:LDPuptbp} holds
under the assumption that both\break $\dE_\rho\psi,\dE_\rho\phi$ are finite.

If, for example, $\dE_\rho\psi$
is infinite, then the above argument can be repeated, replacing $\dE
_\rho\psi$ by a large number $K$, and then letting $K\to\infty$ at the
end. The same reasoning applies to the case where
$\dE_\rho\phi=\infty$. Similarly, if, for example, $b=\infty$ and
$\dE
_\rho\psi>0$, one can replace $b$ above by a large number $K$ and then
let $K\to\infty$ at the end. The same applies to the case $a=\infty$
and $\dE_\rho\phi>0$.
In particular, in all these cases one has that the left-hand side of
\eqref{eq:LDPuptbp} is $-\infty$.
\qed

\emph{Lower bound.}
It is sufficient to prove that for any $\rho\in\cP_s (\cG_*)$ and any
$\delta> 0$,
%
\begin{equation}
\label{eq:lbLDPUG} 
\liminf_{n \to\infty} \frac{ 1}{ n^{1 + \alpha/2} }
\log\dP {{\bigl(\rho_n \in B_\proj( \rho, \delta) \bigr)}}
\geq- b \dE_{\rho} \psi - a \dE_{\rho} \phi.
\end{equation}

In order to prove \eqref{eq:lbLDPUG}, we may assume without loss of
generality that $I (\rho) = b \dE_{\rho} \psi+ a \dE_{\rho} \phi<
\infty$. By monotonicity \eqref{monot}, one has that
\[
\lim_{ j \to\infty} I ( \rho_{\theta_j} ) = I (\rho).
\]
Therefore, since the projective topology is generated from the product
topology on $\prod_{ j \geq1} \cG^{\theta_j}_*$, it is sufficient to
prove \eqref{eq:lbLDPUG} for all
$\rho\in\cP_s (\cG^\theta_*)$, for all $0 < \theta< 1$.
Finally, since the local weak topology is finer than the projective
weak topology, it is enough to prove that for any $0 < \theta< 1$,
$\rho\in\cP_s (\cG^\theta_*)$ and $\delta> 0$,
%
\begin{equation}
\label{toplow} 
\liminf_{n \to\infty} \frac{ 1}{ n^{1 + \alpha/2} }
\log\dP {{\bigl(\rho_n \in B_\loc( \rho, \delta) \bigr)}}
\geq- b \dE_{\rho} \psi - a \dE_{\rho} \phi.
\end{equation}

Let us start with some simple consequences of Assumption~\ref{A1}.
From \eqref{hypoa}, there exists a positive sequence $\eta_n$
converging to $0$ such that, for any $s\geq\veps(n)=1/\log n$,
%
\begin{equation}
\label{toplow0} e^ { - (a + \eta_n ) s^{\alpha} n^{\alpha/2} } \leq\dP\bigl( |X_{12} | \geq s \sqrt n \bigr)
\leq e^ { - ( a - \eta_n)
s^{\alpha} n^{\alpha/2} }.
\end{equation}
In particular, if $s\geq\veps(n)$, then for any $\g>0$, for all $n$
large enough,
\[
\dP \bigl( |X_{12} | \in[s, s + \g) \sqrt n \bigr) \geq
\tfrac{1}2 e^ { - ( a + \eta_n) s^{\alpha} n^{\alpha/2} }.
\]
Therefore, using \eqref{hypoaa}, one finds that there exists a sequence
$a_n\to a$ such that for every $\g>0$, for all $n$ large enough, for
every $z \in\dC$, with $|z| \geq\veps(n)$, $z / |z| \in S_a$,
%
\begin{equation}
\label{eq:prelimlwao} \dP{{\bigl(X_{12}/\sqrt n \in B_{\dC} (z, \g)
\bigr)}} \geq e^ { - a_n
|z|^{\alpha} n^{\alpha/2} },
\end{equation}
where $S_a$ denotes the compact support of the measure $\vartheta_a\in
\cP(\dS^1)$ associated to $X_{12}$, and $B_{\dC} (z,\g)$ is the
Euclidean ball in $\dC$, with center $z$ and radius $\g>0$.

Similarly, there exists a sequence $b_n \to b$
such that for every $\g>0$, for all $n$ large enough, for every $x \in
\dR$, with $|x| \geq\veps(n)$, $x / |x| \in S_b$,
%
\begin{equation}
\label{eq:prelimlwa} \dP{{\bigl(X_{11}/\sqrt n \in B_{\dR} (x, \g)
\bigr)}} \geq e^ { - b_n
|x|^{\alpha} n^{\alpha/2} }.
\end{equation}
%
We remark that \eqref{eq:prelimlwao} and \eqref{eq:prelimlwa} are the
only places where the assumption \eqref{hypoaa} is used in this work.

Since $\r\in\cP_s(\cG_*^\theta)$, there exists a sequence of matrices
$H_n \in\cA_n$,
such that the associated network as in \eqref{ghn} is in $\cG^\theta
$ and
such that $U(H_n)\locweak\r$. In particular,
for $n$ sufficiently large one has
\[
U(H_n) \in B_\loc( \rho, \delta/ 2).
\]
From Lemma~\ref{le:compnorm}, there exists $\g=\g(\d,\theta)>0$ such
that if
$ | \omega_{G_n} (i ) - H_n (i,i )| \leq\g$ and $| \omega_{G_n} (i,j)
- H_n (i,j )| \leq\g$ for all $1\leq i\leq j\leq n$, then $ \r
_n=U(G_n)\in B_\loc( U(H_n), \delta/ 2)$.
Then, by the triangle inequality, for all $n$ large enough,
\begin{eqnarray*}
&&\dP{{\bigl(\rho_n \in B_\loc( \rho, \delta) \bigr)}}\\
&&\qquad
\geq \dP{{\bigl(\rho_n \in B_\loc\bigl( U(H_n)
, \delta/2\bigr) \bigr)}}
\\
&&\qquad \geq\dP \Bigl( \max_{ 1 \leq i \leq n}\bigl | \omega_{G_n} (i ) -
H_n (i,i ) \bigr| \leq\g, \max_{ 1 \leq i < j \leq n} \bigl| \omega
_{G_n} (i,j ) - H_n(i,j ) \bigr| \leq\g \Bigr).
\end{eqnarray*}
Independence of the weights $\omega_{G_n} (i,j )=C_{i,j}$, $1 \leq i
\leq j \leq n$ then gives
\begin{eqnarray*}
&&\dP{{\bigl(\rho_n \in B_\loc( \rho, \delta) \bigr)}} \\
&&\qquad
\geq\prod_{i=1} ^ n \dP{{\bigl({{\bigl|
C_{ii}- H_n (i,i ) \bigr|}} \leq\g\bigr)}} \prod
_{1 \leq i < j \leq n} \dP{{\bigl({{\bigl| C_{ij}- H_n (i,j
) \bigr|}} \leq\g\bigr)}}.
\end{eqnarray*}
Notice that whenever $H_n(i,j)\neq0$ one has $|H_n(i,j)|\geq\theta$,
and thus using \eqref{toplow0} and~\eqref{eq:prelimlwa}
one has for all $i=1,\ldots,n$:
\begin{eqnarray*}
&&\dP{{\bigl({{\bigl| C_{ii}- H_n (i,i ) \bigr|}} \leq\g\bigr)}}
\\
&&\qquad\geq e^ { - b_n n^{\alpha/2} | H_n (i,i) |^\alpha}
 \bigl(\IND{ \bigl(\bigl| H_n(i,i )\bigr | > 0\bigr)} \\
 &&\hspace*{85pt}\qquad{}+
\bigl( 1 - e^{ - c \veps(n)^\alpha n^{\alpha/2}} \bigr)
\IND{ \bigl(\bigl| H_n(i,i )\bigr | = 0
\bigr)} \bigr)
\\
&&\qquad \geq e^ { - b_n n^{\alpha/2} | H_n (i,i) |^\alpha} \bigl( 1 - e^{ -
c \veps(n)^\alpha n^{\alpha/2}} \bigr),
\end{eqnarray*}
where the constant $c$ satisfies $c\geq b/2>0$.
Similarly, using \eqref{eq:prelimlwao}, for all $i\leq j$ and for some
$c\geq a/2>0$:
\[
\dP{{\bigl({{\bigl| C_{ij}- H_n (i,j )\bigr |}} \leq\g\bigr)}}
\geq e^ { - a_n n^{\alpha/2} | H_n (i,j) |^\alpha} \bigl( 1 - e^{ - c \veps(n)^\alpha n^{\alpha/2}} \bigr).
\]
Observe that
\[
\frac{1}n\sum_{1\leq i\leq n} \bigl|H_n (i,i )
\bigr|^\alpha=\dE_{U(H_n)} \psi,\qquad \frac{1}n\sum
_{1\leq i < j\leq n} \bigl|H_n (i,j ) \bigr|^\alpha=
\dE_{U(H_n)} \phi.
\]
Summarizing, using $( 1 - e^{ - c \veps(n)^\alpha n^{\alpha
/2}})^{n^2}\geq1/2$ for $n$ large enough, one finds
%
\begin{equation}
\label{toplow1} \dP{{\bigl(\rho_n \in B_\loc( \rho, \delta)
\bigr)}}  \geq\tfrac{1}2 e^ {
- b_n n^{1 + \alpha/2} \dE_{U(H_n)} \psi} e^ { - a_n n^{1 + \alpha
/2} \dE_{U(H_n)} \phi}.
\end{equation}
Since $\psi$ and $\phi$ are continuous and bounded on $\cG_*^\theta$,
one has $\dE_{U(H_n)} \psi\to\dE_{\r} \psi$ and $\dE_{U(H_n)}
\phi\to
\dE_{\r} \phi$, as $n\to\infty$. Moreover, $a_n\to a$ and $b_n\to b$.
Therefore, \eqref{toplow1}~implies the desired bound \eqref{toplow}.
This concludes the proof of the lower bound.
\end{pf}

The next lemma was used in the proof of the lower bound of Proposition~\ref{prop:LDPnetwork}.
While the estimate is somewhat rough, it is crucial that it is uniform
in the cardinality $n$ of the vertex set.

\begin{lemma}
\label{le:compnorm}
Let $0 < \theta< 1$ and $\delta> 0$.
There exists $\g=\g(\d,\theta) > 0$ such that for any integer $n
\geq
1$, for any networks $G\in\cG$, $H \in\cG_\theta$ with common vertex
set $V = \{1, \ldots, n\}$ such that
%
\begin{equation}
\label{dlocbo1} \max_{ (u, v) \in V^2} \bigl|\omega_G ( u, v) -
\omega_H ( u, v) \bigr|\leq \g,
\end{equation}
then
%
\begin{equation}
\label{dlocbo} \max_{u\in V} \dloc \bigl(\bigl(G(u),u\bigr),
\bigl(H(u),u\bigr) \bigr)\leq\d.
\end{equation}
In particular,
\[
U ( G ) \in B_\loc\bigl( U ( H ), \delta\bigr).
\]
\end{lemma}

\begin{pf}
Each edge of $H$ has a weight bounded above by $\theta^{-1}$. This
implies that in $H$ each path whose total length is bounded by $t>0$,
contains at most $t^2 / \theta^2$ edges. Moreover, $H$ has at most
$\theta^{-4}$ outgoing edges from any vertex. Hence, $H$ has at most $m
= \theta^{-4 t^2 / \theta^2}$ vertices at distance less than $t$ from
any given vertex.
Fix the root $u\in V$ and $t>0$. Therefore, there must exist $ t_0>0$
such that $t /2 < t_0 < t$, and an interval $I = [t_0 - t / (8 m), t_0
+ t / ( 8 m )]$, such that there is no vertex within distance $s\in I$
from $u$ in $H$.

If $e_1, \ldots, e_k$ are the edges on a path in $H$, then provided
that $0 < \g< \theta/2$, one has
\begin{eqnarray*}
&&\Biggl[ \Biggl( \sum_{i=1} ^ k \bigl|
\omega_H ( e_i )\bigr|^{-2} \Biggr)^{1/2}
- \Biggl( \sum_{i=1} ^ k \bigl|
\omega_G ( e_i )\bigr|^{-2} \Biggr)^{1/2}
\Biggr]^2\\
&&\qquad \leq \sum_{i=1}
^ k {{\bigl(\bigl| \omega_H ( e_i )
\bigr|^{-1} - \bigl|\omega_G ( e_i )\bigr|^{-1}
\bigr)}}^2 \leq\frac{4 \g^2 k} { \theta^4}.
\end{eqnarray*}
%
The first inequality follows from the convexity of
$[0,\infty)^2\ni(x,y)\mapsto(\sqrt x - \sqrt y)^2$, which yields
$
((\sum_{i}u_i)^{1/2}-(\sum_{i}v_i)^{1/2} )^2\leq\sum_{i}(u_i^{1/2}-v_i^{1/2})^2$, for any $u,v\in\dR_+^k$. The second
inequality follows from $|\omega_H ( e_i) | \geq\theta$ and the
assumption \eqref{dlocbo1}.
In the worst possible case, one can take $ k = t^2 / \theta^2$ for the
number of edges at distance $t_0$ from $u$.
Together with the previous observation, this shows that if $ 2 \g\sqrt
k / \theta^2\leq t / ( 8m )$, that is, $\g\leq\theta^{3}/(16 m)$,
then the neighborhood of $u$ consisting of vertices within distance
$t_0$ in $G$ and in $H$ have the same vertex set.
From the definition of $\dloc$, this choice of $\g$ in \eqref{dlocbo1}
implies that
\[
\dloc\bigl( \bigl(G(u), u \bigr), \bigl(H(u), u\bigr) \bigr) \leq
\frac{1 }{ 1 + \g^{-1} \wedge
t_0 } \leq\frac{2} t.
\]
Thus, taking $t=2/\d$, one has \eqref{dlocbo}, as soon as, for
example, $\g\leq\theta^{3}/(16 m)=\theta^{3+16/(\d^2\theta^2)}/16$.
From the definition of the L\'evy distance, it immediately follows that
$U ( G ) \in B_\loc{{(U ( H ), \d)}}$.
\end{pf}

%
\begin{remark}
In the proof of Proposition~\ref{prop:LDPnetwork}, we have not appealed
to general results, such as Dawson--G\"artner's theorem, that are
available for projective topologies (see, e.g., \cite{dembo},
Section~4.6). We have, however, crucially used the compactness of $\cP
_s (\cG^*)$ for the projective weak topology. It is not hard to check
that the rate function $I(\rho)$ in \eqref{Irho} is not good for the
weak topology (level sets are not compact).
\end{remark}

\subsection{Spectral measure}

For a network $G = (V,\omega) \in\cG^\theta$, we may define the
bounded linear operator $T$ on the Hilbert space $\ell^2 ( V)$ by
%
\begin{equation}
\label{opT} T e_v = \sum_{ u \in V }
\omega(u,v) e_u
\end{equation}
for any $v \in V$, where $\{e_u, u\in V\}$ denotes the canonical
orthonormal basis of $\ell^2 ( V)$.
$T$ is bounded since
%
\begin{equation}
\label{eq:boundedT} \| T e_v \|^2_2 = \sum
_{u \in V } \bigl|\omega(v,u) \bigr|^2 = \deg(v) \leq
\theta^{-2}.
\end{equation}
Also, since $G$ is Hermitian, $T$ is self-adjoint. We may thus define
the spectral measure at vector $e_v$, see, e.g., \cite{reedsimon}, as
the unique probability measure $\mu_T^{v}$ on $\dR$ such that for any
integer $k \geq1$,
%
\begin{equation}
\label{rootsp} \int x^k \,d\mu_T^{v} = \bigl
\langle e_v, T^k e_v \bigr\rangle.
\end{equation}
Notice that for rooted networks $(G,o)$ with $ G\in\cG^\theta$, then
the associated spectral measure $\mu_T^{o}$ is constant on the
equivalence class of $(G,o)$, so that $\mu_T^{o}$ can be defined as a
measurable map from $\cG_*^\theta$ to $\cP(\dR)$. Thus,
if $\rho\in\cP( \cG^\theta_*)$ for some $\theta>0$, one can define
the spectral measure of $\r$ as
%
\begin{equation}
\label{thetasp} \mu_\rho= \dE_\rho\mu_T
^ {o}.
\end{equation}
%
In particular, consider a Hermitian matrix $H_n\in\cH_n(\dC)$, let
$G_n=G(H_n)$ be the associated network as in \eqref{ghn}, and let
$\rho
_n = U(G_n)$. Then, if $(\psi_1, \ldots, \psi_n)$ is an orthonormal
basis of eigenvectors of $H_n$ with associated eigenvalues $(\lambda_1,
\ldots, \lambda_n)$, by the spectral theorem, for any $v \in\{
1,\ldots,
n \}$,
\[
\mu^{v} _{H_n} = \sum_{i=1}^ n
\bigl|\langle\psi_i, e_v \rangle\bigr|^ 2 \delta
_{\lambda_i},
\]
where $\mu^{v} _{H_n}$ stands for the spectral measure at $v$; see
\eqref{rootsp}. Moreover, the empirical distribution of the eigenvalues
of $H_n$ satisfies
%
\begin{equation}
\label{eq:murhofinite} \mu_{H_n} = \frac{1} n \sum
_{i = 1} ^ n \delta_{\lambda_i} =
\frac{1} n \sum_{v=1} ^ n
\mu^{v} _{H_n} = \mu_{\rho_n}.
\end{equation}
Hence, our definition of spectral measure for a sofic distribution
coincides for finite networks with the empirical distribution of the
eigenvalues.

We turn to the definition of $\mu_\r$ for the case where $\r\in\cP
(\cG
_*)$ but there is no $\theta>0$ such that $\r\in\cP( \cG^\theta_*)$.
In this case, \eqref{thetasp} allows one to define the spectral
measures $\mu_{\rho_\theta}$, where the truncated network $\r
_\theta$
is defined as in Lemma~\ref{le:contproj}.
Next, we shall define the spectral measure $\mu_\r$ as the limit of
$\mu
_{\rho_\theta}$ as $\theta\to0$, provided some extra assumptions are
satisfied.
More precisely, for a rooted network $(G,o)$, $G \in\cG$, and for $\b
>0$, let
%
\begin{equation}
\label{xibeta} \xi_\b(G,o) = \sum_{ v \in V_G}
\bigl| \omega_{G} (o,v)\bigr |^{\b}.
\end{equation}
Since $\xi_\b$ is constant on the equivalence class of $(G,o)$, it can
be seen as a function on $\cG_*$.
For $\b> 0$, $\t>0$, define
\[
\cP_{s,\b,\t} ( \cG_*)=\bigl\{\r\in\cP_{s} ( \cG_*)\dvtx
\dE_{\rho} \xi_\b< \t\bigr\}.
\]
%
Lemma~\ref{le:defspec} below is an 
extension to the weighted case of analogous statements in \cite
{MR2724665,MR2789584}, where spectral measures are defined for random
rooted graphs (with no weights). The first result allows one to define
the spectral measure $\mu_\r$ of any $\rho\in\cP_{s,\b,\t} ( \cG_*)$.

\begin{lemma}
\label{le:defspec}
Let $0 < \b< 2$, $\t>1$ and $\rho\in\cP_{s,\b,\t} ( \cG_*)$. Then
the weak limit
\[
\mu_\rho:= \lim_{\theta\to0} \mu_{\rho_\theta}
\]
exists in $\cP(\dR)$.
\end{lemma}

\begin{pf}
To prove the lemma, we are going to show that the sequence $\mu_{\r
_\theta}$, $\theta\to0$, is Cauchy w.r.t. the metric \eqref{eq:defdist}.

By assumption, there exists a sequence $G_n$ of networks on $\{ 1,
\ldots, n\}$ such that $\r_n \locweak\rho$, where $\rho_n = U(G_n) $.
Call $T_n$ the associated Hermitian matrix. The empirical distribution
of the eigenvalues of $T_n$ satisfies 
$\mu_{T_n} = \mu_{\rho_n}$ by \eqref{eq:murhofinite} applied to $H_n=T_n$.

The truncations $(\rho_n)_\theta$ and $\rho_\theta$ satisfy $(\rho
_n)_\theta\locweak\rho_\theta$ by Lemma~\ref{le:contproj}(ii).
Moreover, for all $\theta>0$,
%
\begin{equation}
\label{eq:convGveps} \mu_{(\rho_n)_\theta} \weak\mu_{\rho_\theta}.
\end{equation}
To prove \eqref{eq:convGveps}, let $T^\theta$ denote the random bounded
self-adjoint operator associated to $\r_\theta$ via \eqref{opT} and let
$T_n^\theta$ be the matrices associated to $(\rho_n)_\theta$. One can
realize these operators on a common Hilbert space $\ell^2(V)$.
Since $(\rho_n)_\theta\locweak\rho_\theta$, from the Skorokhod
representation theorem one can
define a common probability space such that the associated networks
converge locally almost surely, so that a.s. $T_n^\theta e_v\to
T^\theta e_v$, in $\ell^2(V)$, for any $v\in V$.
This implies the strong resolvent convergence; see, for example,
\cite{reedsimon}, Theorem VIII.25(a), and in particular that for any
$v \in
V$, a.s. %
\[
\mu^{v}_{T_n^\theta} \weak\mu^{v}_{T^\theta}.
\]
Then \eqref{eq:convGveps} follows by applying this to $v = o$ and
taking expectation.

Let $T^\theta_n, \wt T^\theta_n $ be the matrices associated to
$(G_n)_\theta$ and $(\wt G_n)_\theta$, respectively, where $(\wt
G_n)_\theta$ is defined according to \eqref{eq:defomthet}, and
$(G_n)_\theta$ according to \eqref{eq:defomthet2}. From \eqref
{eq:boundd}, using the triangle inequality, Lemmas~\ref{le:rank} and~\ref{le:HW},
\[
d ( \mu_{T^\theta_n}, \mu_{T_n} ) \leq\frac{1}n \rank
\bigl( \wt T^\theta _n - T_n\bigr) + \biggl(
\frac{1} n \TR\bigl(\wt T^\theta_n -
T^\theta_n\bigr)^2 \biggr)^{1/2}.
\]
From the definition \eqref{eq:defomthet}, one has
\[
\frac{1}n \rank\bigl( \wt T^\theta_n -
T_n\bigr)  \leq\frac{2 } { n} \sum_{i=1}
^n \IND\bigl( \deg_{G_n} ( i ) \geq \theta^{-2}
-1 \bigr) = 2 \dP_{\rho_n} {{\bigl(\deg_{G} ( o ) \geq\theta
^{-2} -1 \bigr)}}.
\]
From \eqref{eq:defomthet2}, one finds
\begin{eqnarray*}
\frac{1} n \TR\bigl(\wt T^\theta_n -
T^\theta_n\bigr)^2 &\leq&\frac{1} {n} \sum
_{i,j=1}^n \bigl| \omega_{G_n} (i,j)
\bigr|^2 \IND\bigl(\bigl | \omega_{G_n} (i,j) \bigr| \leq2 \theta\bigr)
\IND\bigl( \deg_{G_n} ( i ) \leq \theta^{-2} \bigr)
\\
& = &\dE_{\rho_n} \IND\bigl( \deg_{G} ( o ) \leq
\theta^{-2} \bigr) \sum_{v} \bigl|
\omega_{G} (o,v) \bigr|^2 \IND\bigl( \bigl| \omega_{G}
(o,v) \bigr| \leq2 \theta\bigr).
\end{eqnarray*}
Letting $n$ go to infinity, using $\mu_{T_n^\theta} = \mu_{(\rho
_n)_\theta}$,
and \eqref{eq:convGveps}, one has $d ( \mu_{T^\theta_n}, \mu
_{T_n^{\theta'}} )\to d ( \mu_{\rho_\theta}, \mu_{\rho_{\theta
'}} )$.
Therefore, by
the triangle inequality and the dominated convergence theorem, for any
$0 < \theta' < \theta< 1 / \sqrt2$,
\begin{eqnarray*}
d ( \mu_{\rho_\theta}, \mu_{\rho_{\theta'}} ) & \leq& 4 \dP _{\rho} {{
\bigl(\deg_{G} ( o ) \geq\theta^{-2} /2 \bigr)}}
\\
&&{} + 2 \biggl( \dE_{\rho} \IND\bigl( \deg_{G} ( o ) \leq
\theta^{-2} \bigr) \sum_{v} \bigl| \omega
_{G} (o,v)\bigr |^2 \IND\bigl( \bigl| \omega_{G} (o,v)
\bigr| \leq2 \theta\bigr) \biggr)^{1/2}.
\end{eqnarray*}
Notice that, for $\b\in(0,2)$
%
\begin{equation}
\label{eq:degkappa} \deg_{G} ( o )^{\b/2} = \biggl(\sum
_{v}\bigl | \omega_{G} (o,v)\bigr |^2
\biggr)^{\b/2} \leq\sum_{v}\bigl |
\omega_{G} (o,v) \bigr|^{\b} =\xi_\b( G,o),
\end{equation}
where we use that $\sum_{i=1}^k a_i^r \leq ( \sum_{i=1}^k
a_i
)^r$ for all $a_i\geq0$, $r\geq1$ and $k\in\dN$.
Moreover,
\[
\sum_{v} \bigl| \omega_{G} (o,v)
\bigr|^2 \IND\bigl(\bigl | \omega_{G} (o,v) \bigr| \leq \theta \bigr)
\leq\theta^{ 2 - \b} \xi_\b(G,o).
\]
Hence, from Markov's inequality,
%
\begin{equation}
\label{eq:distcauchy} d ( \mu_{\rho_\theta}, \mu_{\rho_{\theta'}} ) \leq4
\theta^
{\b} \dE _\rho\xi_\b+ 2
\theta^{1 - \b/2} ( \dE_{\rho} \xi_\b )^{1/2}.
\end{equation}
By assumption $\dE_\rho\xi_\b$ is finite. Hence, the sequence $\mu
_{\rho_\theta}$ is Cauchy.
\end{pf}

\begin{lemma}
\label{le:contspec}
For any $\b\in(0,2)$, $\tau> 0$, the map $\rho\mapsto\mu_\rho$
from $ \cP_{s,\b,\tau} ( \cG_*)$ to $\cP(\dR)$ is continuous for the
projective weak topology.
\end{lemma}

\begin{pf}
For any $\theta> 0$, from \eqref{eq:distcauchy},
%
\begin{equation}
\label{eq:distcauchy2} d ( \mu_{\rho_\theta}, \mu_{\rho} ) \leq c \bigl(
\theta^ { \b} + \theta^{1 - \b/2} \bigr),
\end{equation}
with a constant $c = c(\t)>0$.
Hence, from the triangle inequality, if $\rho, \rho' \in\cP_{s,\b
,\tau
} ( \cG_*)$,
\[
d ( \mu_\rho, \mu_{\rho'} ) \leq2 c \bigl(
\theta^ {\b} + \theta^{1 -
\b/2} \bigr) + d ( \mu_{\rho_\theta},
\mu_{\rho'_\theta} ).
\]
Consider a sequence $\r'$ such that $\rho' \projweak\rho$.
If $\rho' \projweak\rho$ then $\rho' _\theta\locweak\rho_\theta$
and, therefore, with the same argument used in the proof of \eqref
{eq:convGveps} above one finds
\[
\mu_{\rho'_\theta} \weak\mu_{\rho_\theta}.
\]
We deduce that
\[
\limsup_{\rho' \projweak\rho} d ( \mu_\rho, \mu_{\rho'} )
\leq 2 c \bigl( \theta^\b+ \theta^{1 - \b/2} \bigr).
\]
Since $\theta> 0$ is arbitrarily small, the statement of the lemma follows.
\end{pf}

\subsection{Large deviations for the empirical spectral measure \texorpdfstring{$\mu_C$}{muC}}
We can apply the previous results to the empirical spectral measure
$\mu
_C$, where $C=C(n)$ is the random matrix defined in \eqref{decompo}.
So far, we have defined $\mu_\r$ for every $\r\in\bigcup_{ 0 < \b
< 2 }
\bigcup_{ \tau> 1} \cP_{s, \b,\tau}( \cG_ * )$.
If $\rho\in\cP_s( \cG_ * )$ but $\r\notin\bigcup_{ 0 < \b< 2 }
\bigcup_{ \tau> 1} \cP_{s, \b,\tau}( \cG_ * )$, then we set
\[
\mu_\rho= \delta_0.
\]

\begin{proposition}\label{ldpmuc}
The empirical spectral measures $\mu_C$ satisfy
an LDP on $\cP(\dR)$ equipped with the weak topology,
with speed $n^{1 + \alpha/2}$ and good rate function $\Phi$ given by
%
\begin{equation}
\label{mucphi} \Phi(\nu)=\inf\bigl\{I (\rho), \r\in\cP_s (\cG_*)
\dvtx \mu_\r=\nu\bigr\},
\end{equation}
where $I(\r)$ is the good rate function in Proposition~\ref{prop:LDPnetwork}.
\end{proposition}

\begin{pf}
Recall that by \eqref{eq:murhofinite}
the network $G_n$ in \eqref{network} satisfies $\r_n=U(G_n)$ and
\[
\mu_{\rho_n} = \mu_C.
\]
Notice that if $c ={{(\frac{a} 2 \wedge b )}}$, then
%
\begin{equation}
\label{muc1} I(\rho) \geq c \dE_\rho\xi_\a,
\end{equation}
where $\xi_\a$ is defined by \eqref{xibeta}.
Hence, by Lemma~\ref{le:contspec}, the map $\rho\mapsto\mu_\rho$ is
continuous on the domain of $I(\rho)$.
We would like to apply a contraction principle to get the LDP for $\mu
_{\rho_n}$ from the LDP for $\rho_n$; see, for example, \cite
{dembo}, Theorem~4.2.1(a). However, a little care is needed here because $\rho
\mapsto\mu_\rho$ is continuous on the set $I(\cdot)<\infty$ only.

We start with the lower bound. Assume that $B$ is an open set in $\cP
(\dR)$.
For each $\tau> 0$, by Lemma~\ref{le:contspec}, the function $f_\tau\dvtx
\rho\mapsto\mu_\rho$ from $\cP_{s,\a,\tau} (\cG_*) \to\cP(\dR
)$ is
continuous. Hence, $f_\tau^{-1} (B) $ is an open subset of $\cP_{s,\a
,\tau} (\cG_*)$, and
\[
\dP( \mu_{\rho_n} \in B )\geq\dP\bigl( \rho_n \in
f_\tau^{-1} (B) \bigr).
\]
From Proposition~\ref{prop:LDPnetwork}, it follows that
\[
- \inf_{ \rho\in\cP_{s,\a,\t} ( \cG_*) \dvtx \mu_\rho\in B} I(\rho) \leq\liminf_{ n \to\infty}
\frac{1} { n^{1 + \alpha/2} } \log \dP( \mu_{\rho_n} \in B ).
\]
Using \eqref{muc1}, one has for some $c>0$:
\[
- \inf_{ \rho\in\cP_{s} ( \cG_*) \dvtx \mu_\rho\in B} I(\rho) \leq ( - c \tau) \vee\liminf
_{ n \to\infty} \frac{1} { n^{1 + \alpha/2}
} \log\dP( \mu_{\rho_n} \in B ).
\]
Letting $\tau$ tend to infinity, we obtain the desired lower bound:
\[
- \inf_{ \nu\in B} \Phi(\nu) \leq\limsup_{ n \to\infty}
\frac{1} {
n^{1 + \alpha/2} } \log\dP( \mu_{\rho_n} \in B ).
\]

To prove the upper bound, assume that $B$ is a closed set in $\cP(\dR
)$. By Lemma~\ref{le:contspec}, $f^{-1}_\tau( B) $ is a closed subset
of $\cP_{s,\a,\t}(\cG_*)$. Write
\[
\dP( \mu_{\rho_n} \in B ) \leq\dP\bigl( \mu_{\rho_n} \in B;
\rho_n \in \cP_{s,\a,\tau} (\cG_*) \bigr) + \dP\bigl(
\rho_n \notin\cP_{s,\a,\tau} (\cG_*) \bigr).
\]
Proposition~\ref{prop:LDPnetwork} yields
\[
\limsup_{ n \to\infty} \frac{1} { n^{1 + \alpha/2} } \log\dP\bigl(
\mu_{\rho_n} \in B; \rho_n \in\cP_{s,\a,\tau} (\cG_*)
\bigr) \leq- \inf_{
\rho\in\cP_{s,\a,\tau} (\cG_*)\dvtx \mu_{\rho} \in B} I(\rho),
\]
and, for some $c>0$:
\[
\limsup_{ n \to\infty} \frac{1} { n^{1 + \alpha/2} } \log\dP\bigl(
\rho_n \notin\cP_{s,\a,\tau} (\cG_*) \bigr) \leq- c \tau.
\]
We have checked that
\[
\limsup_{ n \to\infty} \frac{1} { n^{1 + \alpha/2} } \log\dP( \mu_{\rho_n}
\in B ) \leq- \Bigl[(c \tau)\wedge\inf_{\mu\in B} \Phi(\mu) \Bigr].
\]
Letting $\tau$ tend to infinity, we obtain the desired upper bound.
The fact that $\Phi$ is a good rate function can be seen as in \cite{dembo},
Theorem~4.2.1(a), or, more directly, it follows from Lemma~\ref
{le:goodJ0} below.
\end{pf}
%

\subsection{Proof of Theorem \texorpdfstring{\protect\ref{th:main}}{1.1}}
\label{subsec:main2}
Thanks to Proposition~\ref{prop:exptight}, all we have to show is that
is that the sequence of measures $\mu_{\mathrm{sc}} \boxplus\mu_{C} $ satisfies
a LDP in $\cP(\dR)$ with speed $n^{1+ \alpha/2}$,
with the good rate function $\Phi$ defined in Proposition~\ref{ldpmuc}.
Since the map $\nu\mapsto\mu_{\mathrm{sc}} \boxplus\nu$ is continuous in
$\cP
(\dR)$, the above is an immediate consequence of Proposition~\ref{ldpmuc}
and the standard contraction principle. This completes the proof of
Theorem~\ref{th:main}.

\subsection{On the rate function \texorpdfstring{$\Phi$}{Phi}} \label{subsec:Th23}
We turn to a proof of the properties of the rate function listed in
Theorems \ref{th:rate} and~\ref{th:rateb}.

\begin{lemma}
\label{le:goodJ0}
For any $\b\in(0,2)$,
$\t>1$, for any $\r\in\cP_{s,\b,\t}(\cG_*)$, one has
%
\begin{equation}
\label{bibon} \int|x|^\b\,d \mu_{\rho} (x)\leq
\dE_{\rho} \xi_\b.
\end{equation}
\end{lemma}

\begin{pf}
We use the following Schatten
bound: for all $0 < p \leq2$,
%
\begin{equation}
\label{zhan} \int|x|^p\,d\mu_{A}(x) \leq
\frac{1}{n}\sum_{k=1}^n \Biggl(\sum
_{j=1}^n|A_{kj}|^2
\Biggr)^{{p}/2}
\end{equation}
for every Hermitian matrix $A\in\cH_n(\dC)$.
For a proof, see Zhan \cite{zhan}, proof of Theorem~3.32. For $\r\in
\cP
_{s,\b,\t}(\cG_*)$, there exists a sequence of matrices $H_n$ such that
$\r_n=U(H_n) \locweak\rho$. Let $T^\theta_n$ be the Hermitian matrix
associated to $(H_n)_\theta$, the truncated network. From
\eqref{zhan} and \eqref{eq:degkappa}, one has for all $\theta>0$:
\[
\int|x|^\b\,d \mu_{T_n^\theta} (x) \leq \dE_{\rho_n}
\biggl[ \biggl(\theta^{-2} \wedge\sum_{v}
\bigl| \omega (o,v) \bigr| ^2 \biggr)^{\b/2} \biggr] \leq
\dE_{\rho_n} \bigl( \theta^{-\beta} \wedge\xi_\b( G,o)
\bigr).
\]
For $\theta>0$, the spectral measures $\mu_{T_n^\theta}=\mu_{(\r
_n)_\theta}$ have compact support uniformly in~$n$.
Thus, letting $n$ go to infinity, from \eqref{eq:convGveps} one has
%
\begin{equation}
\label{zhan11} \int|x|^\b\,d \mu_{\rho_\theta} (x) \leq
\dE_{\rho} \xi_\b.
\end{equation}
On the other hand, by definition of $\mu_\r$ (see Lemma~\ref
{le:defspec}), one has $\mu_{\r_\theta}\weak\mu_\r$, $\theta\to
0$ and,
therefore,
\[
\int|x|^\b\,d \mu_{\rho} (x)\leq\liminf
_{\theta\to0}\int|x|^\b\,d \mu _{\rho_\theta} (x).
\]
This proves the claim \eqref{bibon}.
\end{pf}

\begin{pf*}{Proof of Theorem~\ref{th:rate}(a)}
The proof is an immediate consequence of Lemma~\ref{le:goodJ0}. Indeed,
from \eqref{muc1} and the definition of $\Phi$, it suffices to show
that for any $\t>1$, for any $\r\in\cP_{s,\a,\t}(\cG_*)$, one has
%
\begin{equation}
\label{bibone} \int|x|^\alpha\,d \mu_{\rho} (x)\leq
\dE_{\rho} \xi_\a.
\end{equation}
This is the case $\a=\b$ in \eqref{bibon}.
\end{pf*}

\begin{pf*}{Proof of Theorem~\ref{th:rate}(b)}
For $x\in\dR$, let $\frg_x\in\cG_*$ denote the network consisting
of a
single vertex $o$ with weight $\o(o,o)=x$.
If $\nu\in\cP(\dR)$, let $\r\in\cP(\cG_*)$ denote the law $\r
=\int_{\dR
} \d_{\frg_x}\,d\nu(x)$. Notice that
\[
\dE_\r\xi_\a= \int_{\dR}
|x|^\a\,d\nu(x) = m_\a(\nu).
\]
Thus, we can assume $\dE_\r\xi_\a<\infty$, otherwise there is nothing
to prove.
Since we assume $\supp(\vartheta_b)=\{-1,+1\}$, one has that $\r$ is
admissible sofic; see Example~\ref{sinet}, and
$\r\in\cP_{s,\a,\t}(\cG_*)$ for some $\t>1$. The spectral
measure $\mu
_\r$ of $\r$, defined as in Lemma~\ref{le:defspec} is easily seen
to be $\mu_\r=\nu$. Then $\Phi(\nu)\leq I(\r)=b \dE_\r\xi_\a
=b m_\a
(\nu)$.
\end{pf*}

\begin{pf*}{Proof of Theorem~\ref{th:rate}(c)}
Thanks to parts (a) and~(b), all we need to prove is that
%
\begin{equation}
\label{dine1} \Phi(\nu)\leq\frac{a}2 m_\a(\nu),
\end{equation}
for all symmetric probabilities $\nu$ on $\dR$.

For $z\in\dC$, let $\hat{\frg}_z\in\cG_*$ denote the equivalence class
of the two vertex network $(V,\o,o)$, with $V=\{o,1\}$,
$\o(o,1)=z$, $\o(1,o)=\bar z$ and $\o(o,o)=\o(1,1)=0$.
Fix some $e^{i\varphi} \in S_a=\supp(\vartheta_a)$, let
$T$ be a nonnegative random variable with some distribution $\mu_+$ on
$[0,\infty)$, and let $\mu\in\cP(\dC)$ denote the law of
$Te^{i\varphi}$.
The law
\[
\r=\frac{1}2\int_{\dC} (\d_{\hat{\frg}_z}+
\d_{\hat{\frg
}_{\bar z}} )\,d\mu(z),
\]
is sofic; see Example~\ref{dinet}. A simple computation shows that the
spectral measure of $\r$ satisfies
$\mu_\r=\mu_\sym$, where
$\mu_\sym$ denotes the symmetric probability on $\dR$ such that
\[
\int_\dR f(x)\,d\mu_\sym(x)=\frac{1}2
\int_{0}^\infty \bigl(f(x)+f(-x) \bigr) \,d\mu_+(x)
\]
for all bounded measurable $f$.

To prove \eqref{dine1}, let $\nu\in\cP_\sym(\dR)$ and write $\mu
_+$ for
the law of $|X|$ when $X$ has law~$\nu$. Then $\nu=\mu_\sym$ and the
associated $\r$ satisfies $\mu_\r=\nu$. Therefore,
\[
\Phi(\nu)\leq I(\r)= \frac{a}2\int_{0}^\infty
x^\a \,d\mu_+(x) =\frac
{a}2 m_\a(\nu).
\]
\upqed\end{pf*}

\begin{pf*}{Proof of Theorem~\ref{th:rateb}(\textup{a})}
We proceed as in the proof of Theorem~\ref{th:rateb}(b). Here, $S_b=\{
+1\}$, and thus the law $\r=\int_{\dR} \d_{\frg_x}\,d\nu(x)$ that
we used
there is not necessarily admissible sofic. However, it is so if one
assumes $\supp(\nu)\subset\dR_+$. The rest of the argument applies
with no modifications.
\end{pf*}

For the remaining statements, we use the following observation.

\begin{lemma}
\label{integle}
If $\rho\in\cP_{s,\b,\t}(\cG_*)$ for some $\b\in(1,2)$, $\t
>1$, then
%
\begin{equation}
\label{integle1} \int_\dR x \,d\mu_\r(x)=
\dE_\r\o_G(o).
\end{equation}
\end{lemma}

\begin{pf}
By definition of the spectral measure $\mu_{\r_\theta}$[see \eqref
{rootsp}], for every $\theta>0$ one has
\[
\int_\dR x \,d\mu_{\r_\theta}(x) = \dE_{\r_\theta}
\o_G(o) = \dE _{\r} \o _{G_\theta} (o),
\]
where $G_\theta$ is the truncation of $G$; see \eqref{eq:defomthet2}.
The weights $\o_{G_\theta}(o)$ satisfy
$|\omega_{ G_\theta} ( o)| \leq|\omega_{ G} ( o)|$ and, since $\b>1$,
$\dE_{\rho} |\omega_{ G} ( o)| \leq (\dE_\r\xi_\b
)^{1/\b}<\t
^{1/\b}$.
Thus, by the dominated convergence theorem,
\[
\lim_{\theta\to0}\int_\dR x\,d
\mu_{\r_\theta}(x) = \dE_\r\o_G(o).
\]
From \eqref{zhan11}, and the fact that $\b>1$, we know that the
identity map
$x\mapsto x$ is uniformly integrable for $(\mu_{\r_\theta})_{\theta
>0}$. Therefore, by definition of $\mu_\r$ (see Lemma~\ref
{le:defspec}), the limit above also equals $\int_\dR x\,d\mu_\r(x)$.
\end{pf}

\begin{pf*}{Proof of Theorem~\ref{th:rateb}(\textup{b})}
In view of the bound \eqref{dine1}, it suffices to show that if $\r
\in
\cP_s(\cG_*)$ with $\mu_\r=\nu$, then
%
\begin{equation}
\label{cicla} \frac{a} 2 \int|x|^\alpha\,d \mu_{\rho}
(x) \leq I(\rho).
\end{equation}
Thanks to \eqref{muc1}, one may assume that $\r\in\cP_{s,\a,\t
}(\cG_*)$
for some $\t>1$. Moreover, by \eqref{muc1} and \eqref{bibone}, we know
that \eqref{cicla} holds if $b\geq a/2$. If $b<a/2$, we proceed as
follows. Since $\a>1$ here, we may apply Lemma~\ref{integle}, and
obtain that
\[
0=\int_\dR x \,d\nu(x)=\dE_\r
\o_G(o),
\]
where we use the symmetry assumption on $\nu$. Since $S_b=\{+1\}$, one
has that $\o_G(o)\geq0$ and, therefore, $\o_G(o)=0$ $\r$-a.s.
In conclusion, $I(\r)=a\dE_\r\phi=\frac{a}2\dE_\r\xi_\a$, and
the claim
\eqref{cicla} follows from \eqref{bibon}.
\end{pf*}

\begin{pf*}{Proof of Theorem~\ref{th:rateb}(\textup{c})}
Suppose that $I(\r)<\infty$. Then by \eqref{muc1}, one has $\r\in
\cP
_{s,\a,\t}(\cG_*)$ for some $\t>1$. Since $\a>1$, Lemma~\ref
{integle} yields
$\int_\dR x \,d\nu(x)=\dE_\r\o_G(o)$ which, together with the
assumption $\int_\dR x \,d\nu(x)<0$, implies
\[
\dE_\r\o_G(o)<0.
\]
However, $S_b=\{+1\}$ implies that $\dE_\r\o_G(o)\geq0$, a contradiction.
Thus, $I(\r)=+\infty$, for all $\r\in\cP_s(\cG_*)$ such that $\mu
_\r=\nu$.
\end{pf*}

\begin{appendix}
\section{Uniform asymptotic freeness}
\label{sec:UAF}
\subsection{Proof of Theorem \texorpdfstring{\protect\ref{th:UAF}}{2.6}}
Recall the definition \eqref{gstie} of the function $g_\mu\dvtx  \dC
_+\mapsto
\dC_+$, for a given $\mu\in\cP(\dR)$.
Theorem~\ref{th:UAF} is a consequence of the following result.

\begin{theorem}[(Uniform bound in subordination formula)]\label{th:pastur}
Let $Y = (Y_{ij} )_{1 \leq i, j \leq n} \in\cH_n (\dC)$ be a Wigner
random matrix with $\VAR( Y_{12} ) = 1$, $\dE|Y_{12} |^3 < \infty$
and $\dE|Y_{11} |^2< \infty$. There exists a universal constant $c
>0$, such that for any integer $n \geq1$, any $M \in\cH_n (\dC)$, any
$z \in\dC_+$, $\Im( z ) \geq1$,
\[
{{\bigl| \overline g (z) - g_{\mu_M} {{\bigl(z + \overline g (z) \bigr)}} \bigr|}}
\leq c \frac{ {{(\dE|Y_{11} |^2 )}}^{1/2}+ \dE|Y_{12} |^3 }{ n^{1/2} },
\]
where $ \overline g (z) = \dE g_{\mu_{Y / \sqrt n + M}} (z)$.
\end{theorem}
%
Theorem~\ref{th:pastur} is a small generalization of Pastur and
Shcherbina \cite{MR2808038}, Theorem~18.3.1: the main difference here
is that we do not assume that the real and imaginary parts of $Y_{ij}$
are independent. We also allow the mean of the entries to be nonzero.
Note that the rate $1/\sqrt n $ in Theorem~\ref{th:pastur} is not
necessarily optimal with stronger assumptions; see, for example,
\cite{MR2835253}, equation (3.8). We postpone the
proof of Theorem~\ref{th:pastur} to the next subsection. We first check
that it implies Theorem~\ref{th:UAF}. This is done by a simple
contraction argument. For $z \in\dC_+$, we define the $\dC_+ \to\dC
_+$ map,
%
\begin{equation}
\label{eq:defphi} \phi_z \dvtx h \mapsto g_{\mu_M} ( z + h ).
\end{equation}
It is Lipschitz with constant $1/ \Im(z) ^2$. In particular, if $\Im
(z) \geq2$, $\phi_z$ is a contraction with Lipschitz constant $1/ 4$.
Now, it is well known that if $\mu= \mu_M \boxplus\mu_{\mathrm{sc}}$, we have
for all $z \in\dC_+$ the subordination formula,
\[
g_\mu(z) = g_{\mu_M} {{\bigl(z + g_{\mu} (z)
\bigr)}} = \phi_z \bigl( g_{\mu} (z)\bigr),
\]
see Biane \cite{MR1488333}. In particular, if for some probability
measure $\nu\in\cP( \dR)$ and $\veps\geq0$,
\[
{{\bigl| g_\nu(z) - g_{\mu_M} {{\bigl(z + g_{\nu} (z)
\bigr)}} \bigr|}} \leq\veps,
\]
then
\[
{{| g_\mu( z) - g_\nu( z) \bigr|}} \leq\veps+ {{\bigl|
\phi_z \bigl(g_\mu( z)\bigr) - \phi_z
\bigl(g_\nu( z)\bigr) \bigr|}} \leq\veps+ \frac{1}{\Im(z)^2} {{\bigl|
g_\mu( z) - g_\nu( z) \bigr|}}.
\]
So that, if $\Im( z) \geq2$,
\[
{{\bigl| g_\mu( z) - g_\nu( z) \bigr|}} \leq\tfrac{4}{3}
\veps.
\]
Hence, from the definition of the distance $d( \mu, \nu)$ in \eqref
{eq:defdist}, we see that Theorem~\ref{th:UAF} is a corollary of
Theorem~\ref{th:pastur}.

\subsection{Proof of Theorem \texorpdfstring{\protect\ref{th:pastur}}{A.1}: The Gaussian case}
\label{subsec:pasturgauss}

In this subsection, we assume that:
\begin{longlist}[(1)]
\item[(1)]
$G = (\Re(Y_{12}), \Im(Y_{12} )) $ is a centered Gaussian vector in
$\dR^2$ with covariance $K \in\cH_2(\dR)$, $\TR(K) = 1$.
\item[(2)]
$Y_{11} $ is a centered Gaussian in $\dR$ with variance $1$.
\end{longlist}

The proof is a variant of Pastur and Shcherbina \cite{MR2808038}, Lemma~2.2.3 (the main difference is that in \cite{MR2808038}, Lemma~2.2.3, the covariance matrix $K$ is diagonal). We first
recall the Gaussian integration by part formula (see, e.g., \cite
{MR2808038}): for any continuously differentiable function $F\dvtx \dR^2
\mapsto\dR$, with $\dE\|\grad F ( G) \|_2 < \infty$,
%
\begin{equation}
\label{eq:IPPG} \dE F ( G) G = K \dE\grad F ( G).
\end{equation}
We identify $\cH_n ( \dC)$ with $\dR^{n^2}$. Then, if $\Phi\dvtx  \cH
_n( \dC
) \mapsto\dC$ is a continuously differentiable function, we define
$D_{jk} \Phi(X) $ as the derivative with respect to $\Re(X_{jk} )$,
and for $1 \leq j \ne k \leq n$, $D'_{jk} \Phi(X) $ as the derivative
with respect to $\Im(X_{jk} )$.

Define the resolvent $R (X) = ( X - z ) ^{-1}$, $z\in\dC_+$. From the
resolvent formula,
%
\begin{equation}
\label{resolvent} R(X+A)-R(X)=-R(X+A)AR(X),
\end{equation}
valid for any matrix $A\in\cH_n(\dC)$, a standard computation shows that
if $1 \leq j, k \leq n$, and $1 \leq a \ne b \leq n$, then
\[
D_{ab} R_{jk} = - ( R_{ja} R_{bk} +
R_{jb}R_{ak}) \quad\mbox{and}\quad D'_{ab}
R_{jk} = - i ( R_{ja} R_{bk} -
R_{jb}R_{ak}),
\]
while if $1 \leq a \leq n$, then
\[
D_{aa} R_{jk} = - R_{ja} R_{ak}.
\]
Set $X= Y / \sqrt n + M $, so that
\[
R = ( Y / \sqrt n + M - z ) ^{-1}.
\]
Using \eqref{eq:IPPG} we get, for $0 \leq a \ne b \leq n$, and all $j,k$:
%
\begin{eqnarray}
\label{eq:calcul} \dE R_{jk} Y_{ab}& =& \frac{1}{\sqrt n}\dE
\bigl[ K_{11} D_{ab} R_{jk} + K_{12}
D'_{ab} R_{jk} + i K_{21}
D_{ab} R_{jk} + i K_{22} D'_{ab}
R_{jk} \bigr]
\nonumber
\\
& =& - \frac{1}{\sqrt n} \dE \bigl[ (K_{11} - K_{22} + i
K_{12} + i K_{21} ) R_{ja} R_{bk}
\nonumber
\\[-8pt]
\\[-8pt]
\nonumber
&&\hspace*{37pt}{} +
(K_{11} + K_{22} - i K_{12} + i
K_{21} ) R_{jb}R_{ak} \bigr]
\nonumber
\\
& =& - \frac{1}{\sqrt n}\dE{{(\gamma R_{ja} R_{bk} +
R_{jb}R_{ak} )}},\nonumber
\end{eqnarray}
where at the last line, we have used the symmetry of $K$ and $\TR( K )
= 1$, together with the notation
\[
\gamma= K_{11} - K_{22} + 2 i K_{12} = \dE
Y_{ab}^2.
\]
Notice that $|\gamma| \leq1$. Similarly, for $a=b$ one has
%
\begin{equation}
\dE R_{jk}Y_{aa}  = - \frac{1}{\sqrt n}\dE
R_{ja} R_{ak}.\label{eq:calcul2}
\end{equation}
Next, set
\[
G ( z ) = ( M - z ) ^{-1}.
\]
Notice that in this case the dependency of $G(z)$ on $z$ is explicit in
our notation. From the resolvent formula \eqref{resolvent},
\[
R = G(z) - \frac{1}{\sqrt n } R Y G (z).
\]
Hence, for $1 \leq j, k \leq n$, using \eqref{eq:calcul}--\eqref{eq:calcul2},
\begin{eqnarray*}
\dE R_{jk} & =& G(z)_{jk} - \frac{1}{\sqrt n } \sum
_{1 \leq a, b \leq
n} \dE[ R_{j a} Y_{a b}] G
(z)_{b k}
\\
& =& G(z)_{jk} + \frac{\gamma}{ n } \sum_{1 \leq a \ne b \leq n}
\dE [R_{j a} R_{ba}] G (z)_{b k} +
\frac{1}{ n } \sum_{1 \leq a, b \leq n} \dE[R_{jb}R_{aa}
] G (z)_{b
k}.
\end{eqnarray*}
We set
\[
g = g_{\mu_{Y / \sqrt n + M}} (z) = \frac{1} n \sum_{a = 1}^n
R_{aa},\qquad \overline g = \dE g,\qquad \underline g = g - \dE g,
\]
and consider the diagonal matrix $D$ with $D_{jk} = \IND_{ j = k }
R_{jk}$. We find
\[
\dE R = G(z) + \dE[g R] G (z) + \frac{\gamma}{ n } \dE\bigl[R \bigl(
R^\top- D \bigr)\bigr] G (z).
\]
Multiplying on the right-hand side by $G(z) ^{-1} = M - z$ and
subtracting $\overline g R$, one has
\[
\dE R ( M - z - \overline g ) = I + \dE\underline g R + \frac{\gamma}{
n } \dE R
\bigl( R^\top- D \bigr).
\]
Multiplying on the right-hand side by $G(z + \overline g )$,
\[
\dE R = G ( z + \overline g ) + \dE\underline g R G(z + \overline g ) +
\frac{\gamma}{ n } \dE R \bigl( R^\top- D \bigr) G(z + \overline g ).
\]
Finally, multiplying by $\frac{1}n$ and taking the trace,
\[
\overline g  = g_{\mu_M} ( z + \overline g ) + \frac{1}{n}\dE
\underline g \TR\bigl[R G(z + \overline g )\bigr] + \frac{\gamma}{ n^2 } \dE \TR
\bigl[R \bigl( R^\top- D \bigr) G(z + \overline g )\bigr].
\]
As a function of the entries of $Y$, $g$ has Lipschitz constant $ O (
n^{-1} \Im(z )^{-2})$.
This fact can be seen, for example, as in \cite{AGZ}, Lemma~2.3.1.
Since the entries of $Y$ satisfy a Poincar\'e inequality, a standard
concentration bound \cite{MR1849347} implies
\[
\dE|\underline g | = O \bigl( n^{-1} \Im(z )^{-2} \bigr).
\]
Also, since $|\TR( A B )| \leq n \|A \|\|B \|$, we find
\[
\biggl\llvert \frac{1} n \TR R G(z + \overline g ) \biggr\rrvert \leq\Im(
z) ^{-2} \quad\mbox{and} \quad\bigl\llvert \TR R \bigl( R^\top- D
\bigr) G(z + \overline g ) \bigr\rrvert \leq2 n \Im(z) ^{-3}.
\]
This concludes the proof of Theorem~\ref{th:pastur} in the Gaussian case.

\subsection{Proof of Theorem \texorpdfstring{\protect\ref{th:pastur}}{A.1}: The general case}

Let $\underline Y_{ij} = Y_{ij} -\dE Y_{12}$. Then $\underline Y - Y$
has rank at most $1$. Hence, by Lemma~\ref{le:rank},
\[
\bigl| g_{\mu_{Y / \sqrt n + M}} (z) - g_{\mu_{\underline Y / \sqrt n +
M}} (z) \bigr| \leq O \bigl( \bigl(n \Im( z)
\bigr)^{-1} \bigr),
\]
where we have used \eqref{eq:boundd} and the fact that $f ( x ) = ( x -
z) ^{-1}$ has a bounded variation norm of order $\Im( z) ^{-1}$. Also,
we recall that the map $\phi_z$ defined by \eqref{eq:defphi} is Lipschitz
with constant $1/ \Im(z) ^2$. Hence, in order to prove Theorem~\ref
{th:pastur}, we assume without loss of generality that the off-diagonal
entries of the matrix are centered: $\dE Y_{12} = 0$.

We now check that the diagonal entries of $Y$ are negligible. Let $Y'$
be the matrix obtained from $Y$ by setting the diagonal equal to zero:
$Y'_{ij} = \IND_{i \ne j} Y_{ij}$.

\begin{lemma}[(Diagonal entries are negligible)] \label{le:pasturdiag}
For $z \in\dC_+$, $\Im z \geq1$,
\[
\bigl|\dE g_{\mu_{Y / \sqrt n + M}} (z) - \dE g_{\mu_{Y' / \sqrt n + M}} (z) \bigr|= O \bigl( \bigl(
\dE|Y_{11} |^2 / n \bigr)^{1/2} \bigr).
\]
\end{lemma}

\begin{pf}
From \eqref{eq:KRdual}, we find
\begin{eqnarray*}
\bigl|\dE g_{\mu_{Y / \sqrt n + M}} (z) - \dE g_{\mu_{Y' / \sqrt n + M}} (z) \bigr|
 & \leq&
\frac{ \dE W_1 ( \mu_{Y / \sqrt n + M},\mu_{Y' /
\sqrt n + M} ) }{ ( \Im z )^2 }
\\
& \leq&\frac{ \dE W_2 ( \mu_{Y / \sqrt n + M},\mu_{Y' / \sqrt n + M}
) }{ ( \Im z )^2 }.
\end{eqnarray*}
Then by Lemma~\ref{le:HW} using Jensen inequality,
\begin{eqnarray*}
\dE W_2 ( \mu_{Y / \sqrt n + M},\mu_{Y' / \sqrt n + M} ) &\leq&
\frac{1}n \Biggl( \sum_{i=1} ^n
\dE|Y_{ii}|^2 \Biggr)^{1/2} \\
&= &\frac{1}{\sqrt
n}
\bigl(\dE|Y_{11} |^2 \bigr)^{1/2}.
\end{eqnarray*}
\upqed\end{pf}

As a consequence of Lemma~\ref{le:pasturdiag}, we can assume without
loss of generality that the diagonal entries of $Y$ are independent
centered Gaussian with variance~$1$.
By Section~\ref{subsec:pasturgauss}, the conclusion of Theorem~\ref
{th:pastur} holds for the matrix $\widehat Y$ whose off-diagonal
entries are centered Gaussian random variables with covariance is $K$,
where $K$ is the covariance of $Y$, and with diagonal entries centered
Gaussian with variance 1.
Therefore, since the map $\phi_z$ defined by \eqref{eq:defphi} is
Lipschitz, in order to prove Theorem~\ref{th:pastur}, it is sufficient
to establish that
%
\begin{equation}
\label{eq:pasturfin} {{\bigl| \dE g_{\mu_{Y / \sqrt n + M}} (z) -
 \dE g_{\mu_{\widehat Y /
\sqrt n + M}} (z) \bigr|}}
\leq c \frac{ \dE|Y_{12} |^3 }{ n^{1/2} }.
\end{equation}
We may repeat verbatim the interpolation trick in Pastur and Shcherbina
\cite{MR2808038}, Theorem~18.3.1.
Consider the random matrix $\widehat Y$, independent of $Y$, and for $0
\leq t \leq1$, define the matrix
\[
Y ( t) = \sqrt{t} Y + \sqrt{1 - t} \widehat Y.
\]
Set $ R (t) = ( Y ( t) / \sqrt n + M - z I ) ^{-1}$. Then, using the
resolvent equation \eqref{resolvent}
%
\begin{eqnarray}\label{eq:Itrick}
&& g_{\mu_{Y / \sqrt n + M}} (z) - g_{\mu_{\widehat Y / \sqrt n + M}} (z)\nonumber \\
&&\qquad = \frac{1} n \int
_0 ^ 1 \frac{d}{d t} \TR R (t) \,dt
\nonumber
\\
&&\qquad = - \frac{1}{ n ^{3/2} } \int_0^1\TR R(t)
Y'(t) R (t) \,dt
\\
&&\qquad = - \frac{1}{2n ^{3/2} } \int_0^1 \TR R(t)
\biggl( \frac{Y }{\sqrt
t} - \frac{ \widehat Y} {\sqrt{1 - t} } \biggr)R (t) \,dt
\nonumber\\
\nonumber
&&\qquad = - \frac{1}{2n ^{3/2} } \int_0^1 \biggl[
\TR R^2(t) \frac{Y
}{\sqrt
t} - \TR R^2 (t)
\frac{ \widehat Y} {\sqrt{1 - t} } \biggr] \,dt.
\end{eqnarray}
Next, consider the extension of \eqref{eq:IPPG} to arbitrary centered
random variable $G$ with covariance $K$. Namely, for any twice
continuously differentiable function $F\dvtx \dR^2 \mapsto\dR$, with
$\dE
\|\grad F ( G) \|_2 < \infty$ and $ \sup_{x \in\dR^2} \|\operatorname{Hess} F
(x)\| < \infty$, a Taylor expansion gives
\[
\label{eq:IPPG2} \dE F ( G)G = K \dE\grad F ( G) + O \Bigl(\dE\|G
\|_2^3 \sup_{x \in
\dR
^2} \bigl\|\operatorname{Hess} F (x)\bigr\|
\Bigr).
\]
Since $Y$ and $\widehat Y$ have the same first two moments, we get for
all $t\in[0,1]$
\begin{eqnarray*}
&&\dE\TR R^2(t) \frac{Y }{\sqrt t} - \dE\TR R^2 (t)
\frac{ \widehat Y} {
\sqrt{1 - t} } \\
&&\qquad =\sum_{1 \leq j, k \leq n} \dE
R^2(t)_{k j} \frac
{Y_{jk}}{\sqrt t} - \dE R^2(t)_{k j}
\frac{\widehat Y_{jk}}{\sqrt{1-
t}}
\\
&&\qquad \leq c \frac{\dE|Y_{12}|^3}{n} \sum_{1 \leq j, k \leq n} \sup
_{ X \in\cH
_n (\dC), \veps, \veps'} \bigl|D^{\veps}_{jk} D^{\veps'}_{jk}
\bigl( R(X)^2 \bigr)_{kj} \bigr|,
\end{eqnarray*}
where $c>0$ is a constant, and $D_{jk}^\veps D^{\veps'}_{jk} $ ranges
over $D^2_{jk}, {D'}^2_{jk}$ and $D_{jk} D'_{jk}$. However, it follows
from \eqref{eq:calcul}--\eqref{eq:calcul2} that
\[
\bigl|D^{\veps}_{jk} D^{\veps'}_{jk} \bigl(
R(X)^2 \bigr)_{kj}\bigr |
\]
is a finite linear combination of products of $4$
resolvent entries of the form $\prod_{i=1}^4R(X)_{u_iv_i}$. Since
for any $X \in\cH_n (\dC)$, $|R(X)_{jk}|\leq( \Im z ) ^{-1}$, one has
for some new constant $c >0$ and for all $t\in[0,1]$:
\[
\biggl|\dE\TR R^2(t) \frac{Y }{\sqrt t} - \dE\TR R^2 (t)
\frac{
\widehat Y} {\sqrt{1 - t} } \biggr| \leq c n \frac{\dE|Y_{12}|^3}{ (
\Im z ) ^{4}}.
\]
Plugging this last upper bound in \eqref{eq:Itrick} concludes the proof
\eqref{eq:pasturfin} and of Theorem~\ref{th:pastur}.

\section{}
\label{app}

In this section, we collect some standard facts that are repeatedly
used in the main text.
For probability measures $\mu,\mu'\in\cP(\dR)$, the \emph
{Kolmogorov--Smirnov} (KS) distance is defined by
%
\begin{equation}
\label{eq:defdKS} d_{\mathrm{KS}} \bigl( \mu, \mu'\bigr) = \sup
_{ t \in\dR} \bigl| \mu( - \infty, t ] - \mu'( - \infty, t ]
\bigr|.
\end{equation}
The KS distance is closely related to functions with bounded
variations. More precisely, for $f\dvtx \dR\mapsto\dR$ the bounded
variation norm is defined as
\[
\| f \|_{\mathrm{BV}} = \sup\sum_{ k \in\dZ} \bigl| f (
x_{k+1}) - f (x_k) \bigr|,
\]
where the supremum is over all sequence $(x_k)_{ k \in\dZ}$ with
$x_{n} \leq x_{n+1}$. If $f = \IND( ( - \infty, t ) )$, then $\| f \|
_{ BV} = 1$ while if the derivative of $f$ is in $L^1(\dR)$, we have
$\| f \|_{\mathrm{BV}} = \int| f ' (x)|\,dx$.
The KS distance is also given by the variational formula
%
\begin{equation}
\label{eq:KSdual} d_{\mathrm{KS}} \bigl( \mu, \mu'\bigr) = \sup
\biggl\{ \int f\,d \mu- \int f\,d \mu' \dvtx \| f \|_{\mathrm{BV}}
\leq1 \biggr\}.
\end{equation}
%
[Indeed, the functions $H_t = \IND( ( -\infty, t)), t \in\dR$, are
the extremal points of the convex set of functions $f$ with $\| f \|
_{\mathrm{BV}} \leq1$ and the map
$f \to\int f\,d \mu- \int f\,d \mu'$ is linear].

For $p \geq1$ and $\mu,\mu'\in\cP(\dR)$ such that $\int|x|^p\,d
\mu
(x)$ and $\int|x|^p\,d \mu'(x)$ are finite, their $L^p$-Wasserstein
distance is defined as
%
\begin{equation}
\label{eq:defWp} W_p \bigl( \mu, \mu'\bigr) = \biggl(
\inf_\pi\int_{\dR\times\dR} | x - y|
^p\,d \pi(x,y) \biggr)^{{1}/ p},
\end{equation}
where the infimum is over all coupling $\pi$ of $\mu$ and $\mu'$ (i.e.,
$\pi$ is probability measure on $\dR\times\dR$ whose first marginal
is equal to $\mu$ and second marginal is equal to $\mu'$). H\"older's
inequality implies that for $1 \leq p \leq p'$,
$
W_p \leq W_{p'}$.

For any $p \geq1$, if $W_p (\mu_n, \mu)$ converges to $0$ then $\mu_n
\weak\mu$. This follows, for example, from the
Kantorovich--Rubinstein duality
%
\begin{equation}
\label{eq:KRdual} W_1 \bigl(\mu, \mu'\bigr) = \sup \biggl
\{ \int f\,d \mu- \int f\,d \mu' \dvtx \| f \|_{\mathrm{ Lip}} \leq1
\biggr\},
\end{equation}
where $\| f \|_{\mathrm {Lip}}$ denotes the Lipschitz constant of $f$ (see,
e.g., Dudley \cite{MR1932358}, Theorem~11.8.2). 

The following inequality is a standard consequence of interlacing; see,
for example, \cite{MR2567175}, Theorem A.43.

%
\begin{lemma}[(Rank inequality)]\label{le:rank}
If $A$, $B$ in $\cH_n(\dC)$, then
\[
d_{\mathrm{KS}} (\mu_A, \mu_B) \leq
\frac{1}n \rank( A - B ).
\]
\end{lemma}

Next, we recall a very useful estimate which allows one to bound
eigenvalue differences in terms of matrix entries. For a proof see, for
example, \cite{AGZ}, Lemma~2.1.19.

%
\begin{lemma}[(Hoffman--Wielandt inequality)]\label{le:HW}
If $A$, $B$ in $\cH_n(\dC)$, then
\[
W_2 ( \mu_A, \mu_B) \leq\sqrt{
\frac{1} n \TR \bigl[( A-B)^2 \bigr]}.
\]
\end{lemma}
\end{appendix}

\section*{Acknowledgments}
The idea of studying this model came from a discussion with Manjunath
Krishnapur.
We also thank the anonymous referees for several helpful comments.

%
%

%



\printaddresses

\end{document}